\documentclass{article}
\usepackage{amsmath,amssymb}
\usepackage{cite}

\begin{document}

\begin{center}
{\bf Determinantal representations of the quaternion weighted
Moore-Penrose  inverse and corresponding Cramer's rule.}\end{center}
\begin{center}{\bf Ivan Kyrchei }\end{center}\begin{center} Pidstrygach Institute for Applied Problems of Mechanics and Mathematics of NAS of Ukraine,
 Ukraine, kyrchei@online.ua\end{center}

\begin{abstract}
Weighted  singular value decomposition (WSVD)  and  a representation of the  weighted Moore-Penrose inverse  of a quaternion matrix by WSVD have been derived. Using this representation,  limit and determinantal representations  of the  weighted Moore-Penrose inverse  of a quaternion matrix have been obtained   within the framework of the theory of the noncommutative column-row determinants. By using the obtained analogs of the adjoint matrix, we get the Cramer rules for the  weighted Moore-Penrose
solutions of left and right systems of quaternion linear
equations.
\end{abstract}

\section{Introduction}
\newtheorem{corollary}{Corollary}[section]
\newtheorem{theorem}{Theorem}[section]
\newtheorem{lemma}{Lemma}[section]
\newtheorem{definition}{Definition}[section]
\newtheorem{remark}{Remark}[section]
\newcommand{\rank}{\mathop{\rm rank}\nolimits}
\newtheorem{proposition}{Proposition}[section]

Let ${\rm
{\mathbb{R}}}$ and ${\rm
{\mathbb{C}}}$ be the real and complex number fields, respectively.
Throughout the paper, we denote  the set of all $m\times n$ matrices over the
quaternion algebra
\[{\rm {\mathbb{H}}}=\{a_{0}+a_{1}i+a_{2}j+a_{3}k\,
|\,i^{2}=j^{2}=k^{2}=-1,\, a_{0}, a_{1}, a_{2}, a_{3}\in{\rm
{\mathbb{R}}}\}\]
by ${\rm {\mathbb{H}}}^{m\times n}$, and by ${\rm {\mathbb{H}}}^{m\times n}_{r}$ the set of all $m\times n$ matrices over $\mathbb{H}$ with a rank $r$. Let ${\rm M}\left( {n,{\rm {\mathbb{H}}}} \right)$ be the
ring of $n\times n$ quaternion matrices and $ {\bf I}$ be the identity matrix with the appropriate size. For ${\rm {\bf A}}
 \in {\rm {\mathbb{H}}}^{n\times m}$, we denote by ${\rm {\bf A}}^{ *}$, $\rank {\bf A}$ the conjugate transpose (Hermitian adjoint) matrix and the rank
of ${\rm {\bf A}}$.
 The matrix ${\rm {\bf A}} = \left( {a_{ij}}  \right) \in {\rm
{\mathbb{H}}}^{n\times n}$ is Hermitian if ${\rm {\bf
A}}^{ *}  = {\rm {\bf A}}$.

The  definitions of the generalized inverse matrices can be  extended  to quaternion matrices.

The Moore-Penrose inverse of ${\bf A}\in{\rm {\mathbb{H}}}^{m\times n}$, denoted by ${\bf A}^{\dagger}$, is the unique matrix ${\bf X}\in{\rm {\mathbb{H}}}^{n\times m}$ satisfying the following equations \cite{pen},
 \begin{gather}\label{eq1:MP}  {\rm {\bf A}}{\bf X}
{\rm {\bf A}} = {\rm {\bf A}}; \\
                                 \label{eq2:MP}  {\bf X} {\rm {\bf
A}}{\bf X}  = {\bf X};\\
                                  \left( {\rm {\bf A}}{\bf X} \right)^{ *}  = {\rm
{\bf A}}{\bf X}; \\
                                \left( {{\bf X} {\rm {\bf A}}} \right)^{ *}  ={\bf X} {\rm {\bf A}}. \end{gather}

Let Hermitian positive definite matrices ${\bf M}$ and ${\bf N}$ of order $m$ and $n$, respectively, be given. For
 ${\bf A}\in {\mathbb H}^{m\times n}$, \textbf{the weighted Moore-Penrose inverse} of ${\bf A}$ is the unique
solution ${\bf X}={\bf A}^{+}_{M,N}$ of the matrix equations (\ref{eq1:MP}) and (\ref{eq2:MP}) and the following
equations in ${\bf X}$ \cite{pr}:
\[(3M)\,\,({\bf M}{\rm {\bf A}}{\bf X})^{*}  = {\bf M}{\rm
{\bf A}}{\bf X};\,\,(4N)\,\,({\bf N}{\bf X}{\rm {\bf A}})^{*}  = {\bf N}{\bf X}{\rm
{\bf A}}.\]
In particular, when ${\bf M}={\bf I}_m$ and ${\bf N}={\bf I}_n$, the matrix ${\bf X}$ satisfying the  equations (\ref{eq1:MP}), (\ref{eq2:MP}), (3M), (4N) is the  Moore-Penrose inverse ${\bf A}^{\dag}$.

It is known various representations of the weighted Moore-Penrose. In particular, limit representations have been considered in \cite{wei_rep,serg}.
Determinantal representations of the complex (real) weighted Moore-Penrose have been derived
by full-rank factorization in \cite{st1},  by limit representation in  \cite{liu1} using the method first introduced in \cite{ky_lma1},
and by minors in  \cite{liu2}.
A basic method for finding the Moore-Penrose inverse is based on the  singular value decomposition (SVD). It is available for quaternion matrices, (see, e.g. \cite{zha,ky_math_sci}).
In \cite{ky_math_sci,ky_th}, using SVD  of quaternion matrices,  the limit and determinantal representations of the Moore-Penrose inverse over the quaternion skew field have been obtained within the framework of the theory of the noncommutative column-row determinants that have been introduced in
\cite{ky_lma2}.

The  weighted Moore-Penrose inverse ${\bf A}^{\dag}_{M,N}\in {\rm {\mathbb{C}}}^{m\times n} $ can be explicitly expressed by the weighted
singular value decomposition (WSVD) which at first has been obtained in \cite{loan} by Cholesky factorization. In \cite{galba} WSVD of
real matrices with singular weights has been derived using
weighted orthogonal matrices and weighted pseudoorthogonal matrices.

Song at al. \cite{song1,song2} have studied the weighted  Moore-Penrose inverse over the quaternion skew field and obtained its determinantal representation  within the framework of the theory of the column-row determinants. But  WSVD of
quaternion matrices has not been considered and for obtaining  a determinantal representation there was used auxiliary matrices which different from ${\bf A}$, and weights ${\bf M}$ and ${\bf N}$.

The main goals of the paper are  introducing WSVD of quaternion matrices  and representation of the weighted  Moore-Penrose inverse over the quaternion skew field by WSVD, and then by using this representation, obtaining its limit and determinantal representations.

In this paper we shall adopt the following notation.

Let $\alpha : = \left\{
{\alpha _{1} ,\ldots ,\alpha _{k}} \right\} \subseteq {\left\{
{1,\ldots ,m} \right\}}$ and $\beta : = \left\{ {\beta _{1}
,\ldots ,\beta _{k}} \right\} \subseteq {\left\{ {1,\ldots ,n}
\right\}}$ be subsets of the order $1 \le k \le \min {\left\{
{m,n} \right\}}$. By ${\rm {\bf A}}_{\beta} ^{\alpha} $ denote the
submatrix of ${\rm {\bf A}}$ determined by the rows indexed by
$\alpha$, and the columns indexed by $\beta$. Then, ${\rm {\bf
A}}{\kern 1pt}_{\alpha} ^{\alpha}$ denotes the principal submatrix
determined by the rows and columns indexed by $\alpha$.
 If ${\rm {\bf A}} \in {\rm
M}\left( {n,{\rm {\mathbb{H}}}} \right)$ is Hermitian, then by
${\left| {{\rm {\bf A}}_{\alpha} ^{\alpha} } \right|}$ denote the
corresponding principal minor of $\det {\rm {\bf A}}$.
 For $1 \leq k\leq n$, denote by $\textsl{L}_{ k,
n}: = {\left\{ {\,\alpha :\alpha = \left( {\alpha _{1} ,\ldots
,\alpha _{k}} \right),\,{\kern 1pt} 1 \le \alpha _{1} \le \ldots
\le \alpha _{k} \le n} \right\}}$ the collection of strictly
increasing sequences of $k$ integers chosen from $\left\{
{1,\ldots ,n} \right\}$. For fixed $i \in \alpha $ and $j \in
\beta $, let
\[I_{r,\,m} {\left\{ {i} \right\}}: = {\left\{ {\,\alpha :\alpha \in
L_{r,m} ,i \in \alpha}  \right\}}{\rm ,} \quad J_{r,\,n} {\left\{
{j} \right\}}: = {\left\{ {\,\beta :\beta \in L_{r,n} ,j \in
\beta}  \right\}}.\]

 The paper is organized as follows. We start with some basic concepts and results from the theory of  the row-column determinants and  of Hermitian quaternion matrices  in Section 2.
Weighted  singular value decomposition  and  a representation of the  weighted Moore-Penrose inverse  of quaternion matrices by WSVD have been considered in Subsection 3.1, and its limit representations in Subsection 3.2. In Section 4, we  give the determinantal representations of the  weighted Moore-Penrose inverse. In Subsection 4.1, if the matrices ${\bf N}^{-1}{\bf A}^{*}{\bf M}{\bf A}$ and ${\bf A}{\bf N}^{-1}{\bf A}^{*}{\bf M}$ are Hermitian,   and if they are non-Hermitian in  Subsection 4.2.  In Section 5 we
obtain explicit representation formulas of the weighted Moore-Penrose solutions (analogs of Cramer's rule) of the left and right systems of  linear equations over the quaternion skew field. In Section 5, we give
a numerical example to illustrate the main result.
\section{Preliminaries}
 For a quadratic matrix ${\rm {\bf A}}=(a_{ij}) \in {\rm
M}\left( {n,{\mathbb{H}}} \right)$ can be define $n$ row determinants and $n$ column determinants as follows.

Suppose $S_{n}$ is the symmetric group on the set $I_{n}=\{1,\ldots,n\}$.
\begin{definition}\cite{ky_th}
 The $i$-th row determinant of ${\rm {\bf A}}=(a_{ij}) \in {\rm
M}\left( {n,{\mathbb{H}}} \right)$ is defined  for all $i = \overline{1,n} $
by putting
 \begin{gather*}{\rm{rdet}}_{ i} {\rm {\bf A}} =
{\sum\limits_{\sigma \in S_{n}} {\left( { - 1} \right)^{n - r}{a_{i{\kern
1pt} i_{k_{1}}} } {a_{i_{k_{1}}   i_{k_{1} + 1}}} \ldots } } {a_{i_{k_{1}
+ l_{1}}
 i}}  \ldots  {a_{i_{k_{r}}  i_{k_{r} + 1}}}
\ldots  {a_{i_{k_{r} + l_{r}}  i_{k_{r}} }},\\
\sigma = \left(
{i\,i_{k_{1}}  i_{k_{1} + 1} \ldots i_{k_{1} + l_{1}} } \right)\left(
{i_{k_{2}}  i_{k_{2} + 1} \ldots i_{k_{2} + l_{2}} } \right)\ldots \left(
{i_{k_{r}}  i_{k_{r} + 1} \ldots i_{k_{r} + l_{r}} } \right),\end{gather*}
with
conditions $i_{k_{2}} < i_{k_{3}}  < \ldots < i_{k_{r}}$ and $i_{k_{t}}  <
i_{k_{t} + s} $ for $t = \overline{2,r} $ and $s =\overline{1,l_{t}} $.
\end{definition}
\begin{definition}\cite{ky_th}
The $j$-th column determinant
 of ${\rm {\bf
A}}=(a_{ij}) \in {\rm M}\left( {n,{\mathbb{H}}} \right)$ is defined for
all $j =\overline{1,n} $ by putting
 \begin{gather*}{\rm{cdet}} _{{j}}\, {\rm {\bf A}} =
{{\sum\limits_{\tau \in S_{n}} {\left( { - 1} \right)^{n - r}a_{j_{k_{r}}
j_{k_{r} + l_{r}} } \ldots a_{j_{k_{r} + 1} i_{k_{r}} }  \ldots } }a_{j\,
j_{k_{1} + l_{1}} }  \ldots  a_{ j_{k_{1} + 1} j_{k_{1}} }a_{j_{k_{1}}
j}},\\
\tau =
\left( {j_{k_{r} + l_{r}}  \ldots j_{k_{r} + 1} j_{k_{r}} } \right)\ldots
\left( {j_{k_{2} + l_{2}}  \ldots j_{k_{2} + 1} j_{k_{2}} } \right){\kern
1pt} \left( {j_{k_{1} + l_{1}}  \ldots j_{k_{1} + 1} j_{k_{1} } j}
\right), \end{gather*}
\noindent with conditions, $j_{k_{2}}  < j_{k_{3}}  < \ldots <
j_{k_{r}} $ and $j_{k_{t}}  < j_{k_{t} + s} $ for  $t = \overline{2,r} $
and $s = \overline{1,l_{t}}  $.
\end{definition}

Suppose ${\rm {\bf A}}_{}^{i{\kern 1pt} j} $ denotes the submatrix of
${\rm {\bf A}}$ obtained by deleting both the $i$th row and the $j$th
column. Let ${\rm {\bf a}}_{.j} $ be the $j$-th column and ${\rm {\bf
a}}_{i.} $ be the $i$-th row of ${\rm {\bf A}}$. Suppose ${\rm {\bf
A}}_{.j} \left( {{\rm {\bf b}}} \right)$ denotes the matrix obtained from
${\rm {\bf A}}$ by replacing its $j$-th column with the column ${\rm {\bf
b}}$, and ${\rm {\bf A}}_{i.} \left( {{\rm {\bf b}}} \right)$ denotes the
matrix obtained from ${\rm {\bf A}}$ by replacing its $i$-th row with the
row ${\rm {\bf b}}$.
We  note some properties of column and row determinants of a
quaternion matrix ${\rm {\bf A}} = \left( {a_{ij}} \right)$, where
$i \in I_{n} $, $j \in J_{n} $ and $I_{n} = J_{n} = {\left\{
{1,\ldots ,n} \right\}}$.
\begin{proposition}\label{pr:b_into_brak}  \cite{ky_th}
If $b \in {\mathbb{H}}$, then
 $ {\rm{rdet}}_{ i} {\rm {\bf A}}_{i.} \left( {b
\cdot {\rm {\bf a}}_{i.}}  \right) = b \cdot {\rm{rdet}}_{ i} {\rm
{\bf A}}$ and ${\rm{cdet}} _{{i}}\, {\rm {\bf
A}}_{.i} \left( {{\rm {\bf a}}_{.i} b} \right) = {\rm{cdet}}
_{{i}}\, {\rm {\bf A}}  b$ for all $i = \overline {1,n} $.
\end{proposition}
\begin{proposition} \cite{ky_th}
If for  ${\rm {\bf A}}\in {\rm M}\left( {n,{\mathbb{H}}} \right)$\,
there exists $t \in I_{n} $ such that $a_{tj} = b_{j} + c_{j} $\,
for all $j = \overline {1,n}$, then
\[
\begin{array}{l}
   {\rm{rdet}}_{{i}}\, {\rm {\bf A}} = {\rm{rdet}}_{{i}}\, {\rm {\bf
A}}_{{t{\kern 1pt}.}} \left( {{\rm {\bf b}}} \right) +
{\rm{rdet}}_{{i}}\, {\rm {\bf A}}_{{t{\kern 1pt}.}} \left( {{\rm
{\bf c}}} \right), \\
  {\rm{cdet}} _{{i}}\, {\rm {\bf A}} = {\rm{cdet}} _{{i}}\, {\rm
{\bf A}}_{{t{\kern 1pt}.}} \left( {{\rm {\bf b}}} \right) +
{\rm{cdet}}_{{i}}\, {\rm {\bf A}}_{{t{\kern 1pt}.}} \left( {{\rm
{\bf c}}} \right),
\end{array}
\]
\noindent where ${\rm {\bf b}}=(b_{1},\ldots, b_{n})$, ${\rm {\bf
c}}=(c_{1},\ldots, c_{n})$ and for all ${i = \overline {1,n}}$.
\end{proposition}
\begin{proposition} \cite{ky_th}
If for ${\rm {\bf A}}\in {\rm M}\left( {n,{\mathbb{H}}} \right)$\,
 there exists $t \in J_{n} $ such that $a_{i\,t} = b_{i} + c_{i}$
for all $i = \overline {1,n}$, then
\[
\begin{array}{l}
  {\rm{rdet}}_{{j}}\, {\rm {\bf A}} = {\rm{rdet}}_{{j}}\, {\rm {\bf
A}}_{{\,.\,{\kern 1pt}t}} \left( {{\rm {\bf b}}} \right) +
{\rm{rdet}}_{{j}}\, {\rm {\bf A}}_{{\,.\,{\kern 1pt} t}} \left(
{{\rm
{\bf c}}} \right),\\
  {\rm{cdet}} _{{j}}\, {\rm {\bf A}} = {\rm{cdet}} _{{j}}\, {\rm
{\bf A}}_{{\,.\,{\kern 1pt}t}} \left( {{\rm {\bf b}}} \right) +
{\rm{cdet}} _{{j}} {\rm {\bf A}}_{{\,.\,{\kern 1pt}t}} \left( {{\rm
{\bf c}}} \right),
\end{array}
\]
\noindent where ${\rm {\bf b}}=(b_{1},\ldots, b_{n})^T$, ${\rm
{\bf c}}=(c_{1},\ldots, c_{n})^T$ and for all ${j = \overline
{1,n}}$.
\end{proposition}
\begin{remark}\label{rem:exp_det}Let $
{\rm{rdet}}_{{i}}\, {\rm {\bf A}} = {\sum\limits_{j = 1}^{n}
{{a_{i{\kern 1pt} j} \cdot R_{i{\kern 1pt} j} } }} $  and $ {\rm{cdet}} _{{j}}\, {\rm {\bf A}} = {{\sum\limits_{i = 1}^{n}
{L_{i{\kern 1pt} j} \cdot a_{i{\kern 1pt} j}} }}$ for all $i,j =
\overline {1,n}$, where by $R_{i{\kern 1pt} j}$ and $L_{i{\kern 1pt} j} $ denote  the right and left $ij$-th cofactor of ${\rm
{\bf A}}\in {\rm M}\left( {n, {\mathbb{H}}} \right)$, respectively.
It means that  ${\rm{rdet}}_{{i}}\, {\rm
{\bf A}}$ can be expand by right cofactors
  along  the $i$-th row and ${\rm{cdet}} _{j} {\rm {\bf A}}$ can be expand by left cofactors
 along  the $j$-th column, respectively, for all $i, j = \overline {1,n}$.
\end{remark}
The following theorem has a key value in the theory of the column and row
determinants.
\begin{theorem} \cite{ky_th}\label{theorem:
determinant of hermitian matrix} If ${\rm {\bf A}} = \left( {a_{ij}}
\right) \in {\rm M}\left( {n,{\rm {\mathbb{H}}}} \right)$ is Hermitian,
then ${\rm{rdet}} _{1} {\rm {\bf A}} = \cdots = {\rm{rdet}} _{n} {\rm {\bf
A}} = {\rm{cdet}} _{1} {\rm {\bf A}} = \cdots = {\rm{cdet}} _{n} {\rm {\bf
A}} \in {\rm {\mathbb{R}}}.$
\end{theorem}
Since all  column and  row determinants of a
Hermitian matrix over ${\rm {\mathbb{H}}}$ are equal, we can define the
determinant of a  Hermitian matrix ${\rm {\bf A}}\in {\rm M}\left( {n,{\rm
{\mathbb{H}}}} \right)$. By definition, we put
 \begin{equation}\label{eq:def_det_her}\det {\rm {\bf A}}: = {\rm{rdet}}_{{i}}\,
{\rm {\bf A}} = {\rm{cdet}} _{{i}}\, {\rm {\bf A}}, \end{equation}
 for all $i =\overline{1,n}$.
By using its row and column determinants the determinant of a quaternion Hermitian matrix   has properties similar to a usual determinant of a real matrix. These properties are completely explored
in
 \cite{ky_th}
 and can be summarized in the following
theorems.
\begin{theorem}\label{theorem:row_combin} If the $i$-th row of
a Hermitian matrix ${\rm {\bf A}}\in {\rm M}\left( {n,{\rm
{\mathbb{H}}}} \right)$ is replaced with a left linear combination
of its other rows, i.e. ${\rm {\bf a}}_{i.} = c_{1} {\rm {\bf
a}}_{i_{1} .} + \ldots + c_{k}  {\rm {\bf a}}_{i_{k} .}$, where $
c_{l} \in {{\rm {\mathbb{H}}}}$ for all $ l = \overline{1, k}$ and
$\{i,i_{l}\}\subset I_{n} $, then
\[
 {\rm{rdet}}_{i}\, {\rm {\bf A}}_{i \, .} \left(
{c_{1} {\rm {\bf a}}_{i_{1} .} + \ldots + c_{k} {\rm {\bf
a}}_{i_{k} .}}  \right) = {\rm{cdet}} _{i}\, {\rm {\bf A}}_{i\, .}
\left( {c_{1}
 {\rm {\bf a}}_{i_{1} .} + \ldots + c_{k} {\rm {\bf
a}}_{i_{k} .}}  \right) = 0.
\]
\end{theorem}
\begin{theorem}\label{theorem:colum_combin} If the $j$-th column of
 a Hermitian matrix ${\rm {\bf A}}\in
{\rm M}\left( {n,{\rm {\mathbb{H}}}} \right)$   is replaced with a
right linear combination of its other columns, i.e. ${\rm {\bf
a}}_{.j} = {\rm {\bf a}}_{.j_{1}}   c_{1} + \ldots + {\rm {\bf
a}}_{.j_{k}} c_{k} $, where $c_{l} \in{{\rm {\mathbb{H}}}}$ for
all $ l = \overline{1, k}$ and $\{j,j_{l}\}\subset J_{n}$, then
 \[{\rm{cdet}} _{j}\, {\rm {\bf A}}_{.j}
\left( {{\rm {\bf a}}_{.j_{1}} c_{1} + \ldots + {\rm {\bf
a}}_{.j_{k}}c_{k}} \right) ={\rm{rdet}} _{j} \,{\rm {\bf A}}_{.j}
\left( {{\rm {\bf a}}_{.j_{1}}  c_{1} + \ldots + {\rm {\bf
a}}_{.j_{k}}  c_{k}} \right) = 0.
\]
\end{theorem}
\begin{theorem} \cite{ky_th}\label{kyrc12}
If an arbitrary column  of ${\rm {\bf A}}\in {\bf H}^{m\times n} $
is a right  linear combination of its other columns, or an
arbitrary row of ${\rm {\bf A}}^{ * }$ is a left linear
combination of its others, then $\det {\rm {\bf A}}^{ * }{\rm {\bf
A}} = 0.$
\end{theorem}
Moreover, we have the criterion of nonsingularity of a Hermitian matrix.
\begin{theorem} \cite{ky_th}\label{kyrc13}
The right-linearly independence of columns of ${\rm {\bf A}} \in
{\bf H}^{m\times n}$ or the left-linearly independence of rows of
${\rm {\bf A}}^{ *} $ is the necessary and sufficient condition
for $\det {\rm {\bf A}}^{ *} {\rm {\bf A}} \neq 0.$
\end{theorem}
The following theorem about determinantal representation of an
inverse matrix of Hermitian follows immediately from these
properties.
\begin{theorem}
 \cite{ky_th}\label{kyrc10}
 If a Hermitian matrix ${\rm
{\bf A}} \in {\rm M}\left( {n,{\rm {\mathbb{H}}}} \right)$ is such
that $\det {\rm {\bf A}} \ne 0$, then there exist a unique right
inverse  matrix $(R{\rm {\bf A}})^{ - 1}$ and a unique left
inverse matrix $(L{\rm {\bf A}})^{ - 1}$, and $\left( {R{\rm {\bf
A}}} \right)^{ - 1} = \left( {L{\rm {\bf A}}} \right)^{ - 1} =
:{\rm {\bf A}}^{ - 1}$, which possess the following determinantal
representations:
\begin{equation}\label{eq:det_her_inv_R}
  \left( {R{\rm {\bf A}}} \right)^{ - 1} = {\frac{{1}}{{\det {\rm
{\bf A}}}}}
\begin{pmatrix}
  R_{11} & R_{21} & \cdots & R_{n1} \\
  R_{12} & R_{22} & \cdots & R_{n2} \\
  \cdots & \cdots & \cdots & \cdots \\
  R_{1n} & R_{2n} & \cdots & R_{nn}
\end{pmatrix},\end{equation}\begin{equation}\label{eq:det_her_inv_L}
  \left( {L{\rm {\bf A}}} \right)^{ - 1} = {\frac{{1}}{{\det {\rm
{\bf A}}}}}
\begin{pmatrix}
  L_{11} & L_{21} & \cdots  & L_{n1} \\
  L_{12} & L_{22} & \cdots  & L_{n2} \\
  \cdots  & \cdots  & \cdots  & \cdots  \\
  L_{1n} & L_{2n} & \cdots  & L_{nn}
\end{pmatrix}.
\end{equation}
Here $R_{ij}$,  $L_{ij}$ are right and left $ij$-th cofactors of
${\rm {\bf
 A}}$ respectively for all $i,j =
\overline {1,n}$.
\end{theorem}

Due to the noncommutativity of quaternions, there are two types of
eigenvalues.
\begin{definition}
Let ${\rm {\bf A}} \in {\rm M}\left( {n,{\rm {\mathbb{H}}}}
\right)$. A quaternion $\lambda$ is said to be a right eigenvalue
of ${\rm {\bf A}}$ if ${\rm {\bf A}} \cdot {\rm {\bf x}} = {\rm
{\bf x}} \cdot \lambda $ for some nonzero quaternion column-vector
${\rm {\bf x}}$. Similarly $\lambda$ is a left eigenvalue if ${\rm
{\bf A}} \cdot {\rm {\bf x}} = \lambda \cdot {\rm {\bf x}}$.
\end{definition}
The theory on the left eigenvalues of quaternion matrices has been
investigated, in particular, in \cite{hu, so, wo}. The theory on the
right eigenvalues of quaternion matrices is more developed, for further details one
may refer to   \cite{bak, dra, zha}. From this theory we cite
the following propositions.

\begin{proposition}\cite{zha}
Let ${\rm {\bf A}} \in {\rm M}\left( {n,{\rm {\mathbb{H}}}}
\right)$ is Hermitian. Then ${\rm {\bf A}}$ has exactly $n$ real
right eigenvalues.
\end{proposition}
\begin{definition}
Suppose ${\rm {\bf U}} \in {\rm M}\left( {n,{\rm {\mathbb{H}}}}
\right)$ and ${\rm {\bf U}}^{ *} {\rm {\bf U}} = {\rm {\bf U}}{\rm
{\bf U}}^{ *}  = {\rm {\bf I}}$, then the matrix ${\rm {\bf U}}$
is called unitary.
\end{definition}

\begin{proposition} \cite{zha}
Let ${\rm {\bf A}} \in {\rm M}\left( {n,{\rm {\mathbb{H}}}}
\right)$ be given. Then, ${\rm {\bf A}}$ is Hermitian  if and only
if there are a unitary matrix ${\rm {\bf U}} \in {\rm M}\left(
{n,{\rm {\mathbb{H}}}} \right)$ and a real diagonal matrix ${\rm
{\bf D}} = {\rm diag}\left( {\lambda _{{\kern 1pt} 1} ,\lambda
_{{\kern 1pt} 2} ,\ldots ,\lambda _{{\kern 1pt} n}}  \right)$ such
that ${\rm {\bf A}} = {\rm {\bf U}}{\rm {\bf D}}{\rm {\bf U}}^{
*}$, where $\lambda _{ 1},...,\lambda _{ n} $ are right
eigenvalues of ${\rm {\bf A}}$.
\end{proposition}
Suppose ${\rm {\bf A}} \in {\rm M}\left( {n,{\rm {\mathbb{H}}}}
\right) $ is Hermitian and $\lambda \in {\rm {\mathbb {R}}}$ is
its right eigenvalue, then ${\rm {\bf A}} \cdot {\rm {\bf x}} =
{\rm {\bf x}} \cdot \lambda = \lambda \cdot {\rm {\bf x}}$. This
means that all right eigenvalues of a Hermitian matrix are its
left eigenvalues as well. For real left eigenvalues, $\lambda \in
{\rm {\mathbb {R}}}$, the matrix $\lambda {\rm {\bf I}} - {\rm
{\bf A}}$ is Hermitian.
\begin{definition}
If $t \in {\rm {\mathbb {R}}}$, then for a Hermitian matrix ${\rm
{\bf A}}$ the polynomial $p_{{\rm {\bf A}}}\left( {t} \right) =
\det \left( {t{\rm {\bf I}} - {\rm {\bf A}}} \right)$ is said to
be the characteristic polynomial of ${\rm {\bf A}}$.
\end{definition}
The roots of the
characteristic polynomial of a Hermitian matrix are its real left
eigenvalues, which are its right eigenvalues as well.
\begin{definition}Let  ${\bf A}\in {\mathbb{H}}^{n\times n}$ be a Hermitian matrix and $\pi({\bf A}) = \pi$ be the number of positive eigenvalues of ${\bf A}$, $\nu({\bf A}) = \nu$ be
the number of negative eigenvalues of ${\bf A}$, and $\delta({\bf A}) = \delta$ be the number of zero eigenvalues
of ${\bf A}$. Then the ordered triple $\omega =(\pi, \nu, \delta)$ will be called the inertia of ${\bf A}$. We shall
write $\omega = In\, {\bf A}$.
\end{definition}

\begin{definition} A Hermitian matrix ${\bf A}\in {\mathbb{H}}^{n\times n}$, is called  positive (semi)definite if ${\bf x}^{*}{\bf A}{\bf x} > 0(\geq 0)$ for any nonzero
vector ${\bf x}\in {\mathbb{H}}^{n}$.
\end{definition}
The following properties are equivalent to ${\bf A}$ being positive definite and they  can been expanded obviously to quaternion matrices.
\begin{proposition}
All its eigenvalues are positive.
\end{proposition}
\begin{proposition}\label{pr:lead_min}
Its leading principal minors are all positive.
\end{proposition}
Since all leading principal submatrices of a Hermitian matrix are Hermitian, then we can define leading principal minors as determinants of Hermitian submatrices in terms of Eq.(\ref{eq:def_det_her}).

Every positive definite matrix
${\bf A}\in {\mathbb{H}}^{n\times n}$ has a unique  square root defined by ${\bf A}^{\frac{1}{2}}$. It means, if ${\rm {\bf A}} = {\rm {\bf U}}{\rm {\bf D}}{\rm {\bf U}}^{
*}$ then ${\rm {\bf A}}^{\frac{1}{2}} = {\rm {\bf U}}{\rm {\bf D}}^{\frac{1}{2}}{\rm {\bf U}}^{
*}$.

We have \cite{zha,ky_math_sci} the following theorem about  the singular value decomposition (SVD) of quaternion matrices.
\begin{theorem}
(SVD) Let ${\rm {\bf A}} \in {\rm {\mathbb{H}}}_{r}^{m\times n} $.
There exist unitary quaternion matrices ${\rm {\bf V}} \in {\rm
{\mathbb{H}}}^{m\times m}$ and ${\rm {\bf W}} \in {\rm
{\mathbb{H}}}^{n\times n}$ such that
$
 {\rm {\bf A}} = {\rm {\bf V}}{\rm {\bf \Sigma}
}{\rm {\bf W}}^{ *},
$
where ${\rm {\bf \Sigma} } = \begin{pmatrix}
                   {\bf D}& {\bf 0} \\
                   {\bf 0} & {\bf 0} \\
                 \end{pmatrix} \in {\rm
{\mathbb{H}}}_{r}^{m\times n} $, and ${\rm {\bf D}}_{r} = {\rm diag}\left( {\sigma _{1} ,\sigma _{2}
,\ldots ,\sigma _{r}}  \right), \sigma _{1} \ge \sigma _{2} \ge
\ldots \ge \sigma _{r} > 0$, and $\sigma^{2} _{i} $ is the  nonzero eigenvalues of ${\rm {\bf
A}}^{*}{\rm {\bf
A}}$ for all $i=1,...,r$. Then $
 {\rm {\bf A}}^{\dagger} = {\rm {\bf W}}{\rm {\bf \Sigma}
}^{\dagger}{\rm {\bf V}}^{ *},
$ where ${\rm {\bf \Sigma} } = \begin{pmatrix}
                   {\bf D}^{-1}& {\bf 0} \\
                   {\bf 0} & {\bf 0} \\
                 \end{pmatrix}.$
\end{theorem}
In \cite{ky_th}, using the singular value decomposition  of quaternion matrices,  the limit and determinantal representations of the Moore-Penrose inverse over the quaternion skew field have been obtained as follows.
\begin{lemma}\cite{ky_th}\label{lem:lim_rep_A1}
If ${\rm {\bf A}} \in {\rm {\mathbb{H}}}^{m\times n}$ and ${\rm
{\bf A}}^{\dag}$ is its Moore-Penrose  inverse, then ${\rm {\bf
A}}^{ +}  = {\mathop {\lim} \limits_{\alpha \to 0}} {\rm {\bf
A}}^{ * }\left( {{\rm {\bf A}}{\rm {\bf A}}^{ *}  + \alpha {\rm
{\bf I}}} \right)^{ - 1} = {\mathop {\lim} \limits_{\alpha \to 0}}
\left( {{\rm {\bf A}}^{
* }{\rm {\bf A}} + \alpha {\rm {\bf I}}} \right)^{ - 1}{\rm {\bf
A}}^{ *} $, where $\alpha \in {\rm {\mathbb {R}}}_{ +}  $.
\end{lemma}
\begin{theorem}\cite{ky_th}\label{th:det_rep_A1}
If ${\rm {\bf A}} \in {\rm {\mathbb{H}}}_{r}^{m\times n} $, then
the Moore-Penrose inverse  ${\rm {\bf A}}^{ \dag} = \left( {a_{ij}^{ \dag} } \right) \in
{\rm {\mathbb{H}}}_{}^{n\times m} $ possess the following determinantal representations:
\begin{equation}
\label{eq:det_rer_A_c} a_{ij}^{ \dag}  = {\frac{{{\sum\limits_{\beta \in
J_{r,\,n} {\left\{ {i} \right\}}} {{\rm{cdet}} _{i} \left( {\left(
{{\rm {\bf A}}^{ *} {\rm {\bf A}}} \right)_{\,. \,i} \left( {{\rm
{\bf a}}_{.j}^{ *} }  \right)} \right){\kern 1pt} {\kern 1pt}
_{\beta} ^{\beta} } } }}{{{\sum\limits_{\beta \in J_{r,\,\,n}}
{{\left| {\left( {{\rm {\bf A}}^{ *} {\rm {\bf A}}} \right){\kern
1pt}  _{\beta} ^{\beta} }  \right|}}} }}},
\end{equation}
or
\begin{equation}
\label{eq:det_rer_A_r} a_{ij}^{\dag}  = {\frac{{{\sum\limits_{\alpha \in
I_{r,m} {\left\{ {j} \right\}}} {{\rm{rdet}} _{j} \left( {({\rm
{\bf A}}{\rm {\bf A}}^{ *} )_{j\,.\,} ({\rm {\bf a}}_{i.\,}^{ *}
)} \right)\,_{\alpha} ^{\alpha} } }}}{{{\sum\limits_{\alpha \in
I_{r,\,m}}  {{\left| {\left( {{\rm {\bf A}}{\rm {\bf A}}^{ *} }
\right){\kern 1pt}  _{\alpha} ^{\alpha} } \right|}}} }}}.
\end{equation}
\end{theorem}

\begin{lemma}\cite{ky_th} \label{lem:AAdef} Let ${\bf A}\in {\mathbb{H}}_{r}^{m\times n}$, then ${\bf A}^{*}{\bf A}$ and  ${\bf A}{\bf A}^{*}$ are both positive (semi)definite, and $r$ nonzero eigenvalues of ${\rm {\bf
A}}^{*}{\rm {\bf
A}}$ and ${\rm {\bf
A}}{\rm {\bf
A}}^{*}$ coincide.
\end{lemma}{\bf Proof}. The proof of the second part immediately follows from the singular value decomposition  of ${\bf A}\in {\mathbb{H}}_{r}^{m\times n}$.$\blacksquare$

\begin{lemma}\cite{ky_th}\label{lem:char_her}
If ${\rm {\bf A}} \in {\rm M}\left( {n,{\rm {\mathbb{H}}}}
\right)$ is Hermitian, then $p_{{\rm {\bf A}}}\left( {t} \right) =
t^{n} - d_{1} t^{n - 1} + d_{2} t^{n - 2} - \ldots + \left( { - 1}
\right)^{n}d_{n}$, where $d_{k} $ is the sum of principle minors
of ${\rm {\bf A}}$ of order $k$, $1 \le k < n$, and $d_{n}=\det
{\rm {\bf A}}$.
\end{lemma}

\begin{definition}
A  square matrix  ${\bf Q}\in {\mathbb{H}}^{m\times m}$  is called  ${\bf H}$-weighted unitary  (unitary with weight  ${\bf H}$  )  if  ${\bf Q}^{*}{\bf H}{\bf Q}={\bf I}_{m}$,  where  ${\bf I}_{m}$  is  the  identity matrix.
\end{definition}
The following well-known two facts (see,  e.g. \cite{ho}) on positive definite and Hermitian matrices and their product obviously can be extended to quaternion matrices.
\begin{lemma}\label{lem:ho_AB} Let ${\bf A}\in {\mathbb{H}}^{n\times n}$ be positive definite and  ${\bf B}\in {\mathbb{H}}^{n\times n}$ be Hermitian matrices, respectively. Then ${\bf A}{\bf B}$ is a diagonalizable matrix, it's all eigenvalues are real, and $In{\bf A}{\bf B} = In {\bf A}$.
\end{lemma}
\begin{lemma}\label{lem:ho_diag} Let ${\bf A}\in {\mathbb{H}}^{n\times n}$ be positive definite and  ${\bf B}\in {\mathbb{H}}^{n\times n}$ be Hermitian matrices, respectively. Then there exists nonsingular  ${\bf C}\in {\mathbb{H}}^{n\times n}$ such that ${\bf C}^{*}{\bf A}{\bf C}={\bf I}_{n}$,  and ${\bf C}^{*}{\bf B}{\bf C} = {\bf \Lambda}$, where ${\bf \Lambda}$  is  a  diagonal  matrix.
\end{lemma}

\section{Weighted  singular value decomposition  and
representations of the  weighted Moore-Penrose inverse  of quaternion matrices}

\subsection{Representations of the  weighted Moore-Penrose inverse  of quaternion matrices by WSVD}
Denote ${\bf A}^{\sharp}={\bf N}^{-1}{\bf A}^{*}{\bf M}$. Now, we prove the following theorem about  the weighted  singular value decomposition (WSVD) of quaternion matrices. We give the proof of the theorem that different from analogous for real matrices in \cite{loan}, and this method has more similar  manner to \cite{galba}, where WSVD of  ${\bf A}\in {\mathbb{R}}^{m\times n}$ with positive definite weights  ${\bf B}$ and ${\bf C}$ has been described as $ {\bf A}={\bf U}{\bf D}
{\bf V}^{T}{\bf C}$.
\begin{theorem}\label{theor_weig_A} Let ${\bf A}\in {\mathbb{H}}^{m\times n}_{r}$, and ${\bf M}$ and ${\bf N}$ be positive definite matrices  of order $m$ and $n$, respectively.  Then there exist ${\bf U}\in {\mathbb{H}}^{m\times m}$, ${\bf V}\in {\mathbb{H}}^{n\times n}$ satisfying ${\bf U}^{*}{\bf M}{\bf U}={\bf I}_{m}$ and ${\bf V}^{*}{\bf N}^{-1}{\bf V}={\bf I}_{n}$ such that
\begin{equation}\label{eq:weig_A}
{\bf A}={\bf U}{\bf D}
{\bf V}^{*},
\end{equation}
where ${\bf D}=\left(
                 \begin{array}{cc}
                   {\bf \Sigma} & {\bf 0} \\
                   {\bf 0} & {\bf 0} \\
                 \end{array}
               \right)$,
 ${\bf \Sigma}=diag(\sigma_{1}, \sigma_{2},...,\sigma_{r})$, $\sigma_{1}\geq \sigma_{2}\geq...\geq \sigma_{r}>0$ and $\sigma^{2}_{i}$ is the nonzero eigenvalues of ${\bf A}^{\sharp}{\bf A}$ or  ${\bf A}{\bf A}^{\sharp}$, which coincide.
\end{theorem}
{\bf Proof}. First, consider ${\bf A}^{\sharp}{\bf A}={\bf N}^{-1}{\bf A}^{*}{\bf M}{\bf A}$. Since ${\bf A}^{*}{\bf M}{\bf A}= ({\bf M}^{\frac{1}{2}}{\bf A})^{*}{\bf A}{\bf M}^{\frac{1}{2}}$, then, by Lemma \ref{lem:AAdef}, ${\bf A}^{*}{\bf M}{\bf A}$ is Hermitian positive semidefinite, and by Lemma \ref{lem:ho_AB} all  eigenvalues of ${\bf A}^{\sharp}{\bf A}$ are nonnegative. Denote them by $\sigma_{i}^{2}$, where $\sigma_{1}\geq...\geq \sigma_{n}\geq 0$.

Denote ${\bf L}={\bf A}^{\sharp}{\bf A}$. Since $ {\bf L}{\bf N}^{-1}= {\bf N}^{-1}{\bf A}^{*}{\bf M}{\bf A}{\bf N}^{-1}$ is Hermitian and there exists a nonsingular  ${\bf V}\in {\mathbb{H}}^{n\times n}$ such that ${\bf V}^{*}{\bf N}^{-1}{\bf V}={\bf I}_{n}$, then by
Lemma \ref{lem:ho_diag},
\begin{equation}\label{eq:sym_ln}
{\bf V}^{*}{\bf L}{\bf N}^{-1}{\bf V}={\bf \Lambda},
\end{equation}
where $ {\bf V}$ is unitary with weight  ${\bf N}^{-1}$, and  ${\bf \Lambda}$  is  a  diagonal  matrix.

 It follows from  ${\bf L}={\bf N}^{-1}{\bf V}{\bf \Lambda}{\bf V}^{*}=\left({\bf V}^{*}\right)^{-1}{\bf \Lambda}{\bf V}^{*}$ that ${\bf \Lambda}\equiv {\bf \Sigma}_{1}^{2} $,
  where ${\bf \Sigma}_{1}^{2} $ is diagonal with eigenvalues of ${\bf A}^{\sharp}{\bf A}$ on the principal diagonal,
   $\sigma_{ii}^{2} =\sigma_{ i}^{2}$ for all $i=1,...,n$. Since $\rank\,{\bf L}=\rank\,{\bf A}=r$, then the number of nonzero  diagonal  elements of  ${\bf \Sigma}_{1}^{2}$  is equal $r$. Also, we note that
   \begin{multline}\label{eq:sym_ln_1}
{\bf V}^{*}{\bf L}{\bf N}^{-1}{\bf V}={\bf V}^{*}{\bf N}^{-1}{\bf A}^{*}{\bf M}{\bf A}{\bf N}^{-1}{\bf V}=\\
{\bf V}^{-1}{\bf A}^{*}\left({\bf U}^{*}\right)^{-1}{\bf U}^{-1}{\bf A}\left({\bf V}^{*}\right)^{-1}=
{\bf \Sigma}_{1}^{2}.
\end{multline}

Consider the following matrix,
\begin{equation}\label{eq:P_1}{\bf P}={\bf M}^{\frac{1}{2}}{\bf A} {\bf N}^{-1}{\bf V}\in {\mathbb{H}}^{m\times n}.
\end{equation}
By  virtue of   (\ref{eq:sym_ln}),
\begin{equation}\label{eq:PP}
{\bf P}^{*}{\bf P}=\left({\bf V}^{*} {\bf N}^{-1} {\bf A}^{*}{\bf M}^{\frac{1}{2}} \right)  {\bf M}^{\frac{1}{2}}{\bf A} {\bf N}^{-1}{\bf V}={\bf \Sigma}_{1}^{2}.
\end{equation}

Let us  introduce  the  following  $m\times n$  matrix ${\bf D}\in {\mathbb{H}}^{m\times n}$,
\begin{equation}\label{eq:D}{\bf D}=\begin{pmatrix}
{\bf \Sigma} & {\bf 0} \\
{\bf 0} &{\bf 0}
\end{pmatrix},
\end{equation}
where ${\bf \Sigma} \in {\mathbb{H}}^{r\times r}$ is a diagonal matrix with $\sigma_{1}\geq \sigma_{2}\geq...\geq \sigma_{r}>0$ on the principal diagonal.
Then,
\begin{equation}\label{eq:P_2}{\bf P}={\bf M}^{\frac{1}{2}}{\bf U}{\bf D}.
\end{equation}
By (\ref{eq:P_1}) and (\ref{eq:P_2}), we have $ {\bf M}^{\frac{1}{2}}{\bf A} {\bf N}^{-1}{\bf V}={\bf M}^{\frac{1}{2}}{\bf U}{\bf D}$. Due to the equality $\left({\bf N}^{-1}{\bf V}\right)^{-1}={\bf V}^{*}$, it follows (\ref{eq:weig_A}).

Now we shall prove  (\ref{eq:weig_A}), where  $\sigma^{2}_{i}$ is the nonzero eigenvalues of   ${\bf A}{\bf A}^{\sharp}={\bf A}{\bf N}^{-1}{\bf A}^{*}{\bf M}$. Since ${\bf A}{\bf N}^{-1}{\bf A}^{*}$ and ${\bf M}$ are respectively Hermitian positive semidefinite and definite, then by  by Lemma \ref{lem:ho_AB} all  eigenvalues of ${\bf A}{\bf A}^{\sharp}$ are nonnegative. Primarily, denote them by $\tau_{i}^{2}$, where $\tau_{1}\geq...\geq \tau_{m}\geq 0$, and
denote ${\bf Q}={\bf A}{\bf A}^{\sharp}$. Since $ {\bf M}{\bf Q}= {\bf M}{\bf A}{\bf N}^{-1}{\bf A}^{*}{\bf M}$ is Hermitian and there exists a nonsingular  ${\bf U}\in {\mathbb{H}}^{m\times m}$ such that ${\bf U}^{*}{\bf M}{\bf U}={\bf I}_{m}$, then by
Lemma \ref{lem:ho_diag},
\begin{equation}\label{eq:sym_mq}
{\bf U}^{*}{\bf M}{\bf Q}{\bf U}={\bf \Omega},
\end{equation}
where $ {\bf U}$ is unitary with weight  ${\bf M}$, and  ${\bf \Omega}$  is  a  diagonal  matrix.

 It follows from  ${\bf Q}={\bf U}{\bf \Omega}{\bf U}^{*}{\bf M}={\bf U}{\bf \Omega}{\bf U}^{-1}$ that ${\bf \Omega}\equiv {\bf \Sigma}_{2}^{2} $,
  where ${\bf  \Sigma}_{2}^{2} $ is diagonal with eigenvalues of ${\bf A}{\bf A}^{\sharp}$ on the principal diagonal,
   $\tau_{ii}^{2} =\tau_{ i}^{2}$ for all $i=1,...,m$. Since $\rank\,{\bf Q}=\rank\,{\bf A}=r$, then the number of nonzero  diagonal  elements of  ${\bf  \Sigma}_{2}^{2}$  is equal $r$. Also, we have
 \begin{multline}\label{eq:sym_mq_1}
{\bf U}^{*}{\bf M}{\bf Q}{\bf U}={\bf U}^{*}{\bf M}{\bf A}{\bf N}^{-1}{\bf A}^{*}{\bf M}{\bf U}=\\
{\bf U}^{-1}{\bf A}\left({\bf V}^{*}\right)^{-1}{\bf V}^{-1}{\bf A}^{*}\left({\bf U}^{*}\right)^{-1}=
{\bf  \Sigma}_{2}^{2}.
\end{multline}
 Comparing (\ref{eq:sym_ln_1}) and  (\ref{eq:sym_mq_1}), and due to Lemma \ref{lem:AAdef}, we have that $r$ nonzero eigenvalues of   ${\bf A}{\bf A}^{\sharp}$ coincide with $r$ nonzero eigenvalues of   ${\bf A}^{\sharp}{\bf A}$, i.e. $\sigma_{ i}^{2}=\tau_{ i}^{2}$ for all $i=1,...,r$.

 Consider the following matrix,
\begin{equation}\label{eq:S_1}{\bf S}={\bf U}^{*}{\bf M}{\bf A} {\bf N}^{-\frac{1}{2}}\in {\mathbb{H}}^{m\times n}.
\end{equation}
By  virtue of   (\ref{eq:sym_mq}),
\begin{equation}\label{eq:SS}
{\bf S}{\bf S}^{*}={\bf U}^{*}{\bf M}{\bf A} {\bf N}^{-\frac{1}{2}}\left( {\bf N}^{-\frac{1}{2}}{\bf A}^{*}{\bf M}{\bf U} \right) ={\bf \Sigma}_{2}^{2}.
\end{equation}

 Consider again the matrix ${\bf D}\in {\mathbb{H}}^{m\times n}$ from (\ref{eq:D}).
Then,
\begin{equation}\label{eq:S_2}{\bf S}={\bf D}{\bf V}^{*}{\bf N}^{-\frac{1}{2}}.
\end{equation}
By (\ref{eq:S_1}) and (\ref{eq:S_2}), we have $ {\bf U}^{*}{\bf M}{\bf A} {\bf N}^{-\frac{1}{2}}={\bf D}{\bf V}^{*}{\bf N}^{-\frac{1}{2}}$. From this, due to  $\left({\bf U}^{*}{\bf M}\right)^{-1}={\bf U}$,  we   again obtain (\ref{eq:weig_A}).$\blacksquare$

Now, we prove the following theorem about a representation of ${\bf A}^{\dag}_{M,N}$
by WSVD of ${\bf A}\in {\mathbb{H}}^{m\times n}_{r}$ with weights ${\bf M}$ and ${\bf N}$.
\begin{theorem}\label{tm_weig_A+} Let ${\bf A}\in {\mathbb{H}}^{m\times n}_{r}$,  ${\bf M}$ and ${\bf N}$ be positive definite matrices  of order $m$ and $n$, respectively.  There exist ${\bf U}\in {\mathbb{H}}^{m\times m}$, ${\bf V}\in {\mathbb{H}}^{n\times n}$ satisfying ${\bf U}^{*}{\bf M}{\bf U}={\bf I}_{m}$ and ${\bf V}^{*}{\bf N}^{-1}{\bf V}={\bf I}_{n}$ such that
${\bf A}={\bf U}{\bf D}
{\bf V}^{*}$, where ${\bf D}=\left(
                 \begin{array}{cc}
                   {\bf \Sigma} & {\bf 0} \\
                   {\bf 0} & {\bf 0} \\
                 \end{array}
               \right)$.
Then the weighted Moore-Penrose inverse ${\bf A}^{\dag}_{M,N}$ can be represented
\begin{equation}\label{eq:wsvd_A+}
{\bf A}^{\dag}_{M,N}={\bf N}^{-1}{\bf V}\left(
                 \begin{array}{cc}
                   {\bf \Sigma}^{-1} & {\bf 0} \\
                   {\bf 0} & {\bf 0} \\
                 \end{array}
               \right)
{\bf U}^{*}{\bf M},
\end{equation}
where ${\bf \Sigma}=diag(\sigma_{1}, \sigma_{2},...,\sigma_{r})$, $\sigma_{1}\geq \sigma_{2}\geq...\geq \sigma_{r}>0$ and $\sigma^{2}_{i}$ is the nonzero eigenvalues of ${\bf A}^{\sharp}{\bf A}$ or  ${\bf A}{\bf A}^{\sharp}$, which coincide.
\end{theorem}
{\bf Proof}. To prove the theorem we shall show that ${\bf X}={\bf A}^{\dag}_{M,N}$ expressed by (\ref{eq:wsvd_A+}) satisfies the equations (\ref{eq1:MP}), (\ref{eq2:MP}), (3N), and (4M).
\[
1)\,\, {\rm {\bf A}}{\bf X}
{\rm {\bf A}} = {\bf U}\left(
                 \begin{array}{cc}
                   {\bf \Sigma} & {\bf 0} \\
                   {\bf 0} & {\bf 0} \\
                 \end{array}
               \right)
{\bf V}^{*}{\bf N}^{-1}{\bf V}\left(
                 \begin{array}{cc}
                   {\bf \Sigma}^{-1} & {\bf 0} \\
                   {\bf 0} & {\bf 0} \\
                 \end{array}
               \right)
{\bf U}^{*}{\bf M}{\bf U}\left(
                 \begin{array}{cc}
                   {\bf \Sigma} & {\bf 0} \\
                   {\bf 0} & {\bf 0} \\
                 \end{array}
               \right)
{\bf V}^{*}={\rm {\bf A}},\]
\begin{multline*}
2)\,\,{\bf X} {\rm {\bf A}}
{\bf X} ={\bf N}^{-1}{\bf V}\left(
                 \begin{array}{cc}
                   {\bf \Sigma}^{-1} & {\bf 0} \\
                   {\bf 0} & {\bf 0} \\
                 \end{array}
               \right)
{\bf U}^{*}{\bf M} {\bf U}\left(
                 \begin{array}{cc}
                   {\bf \Sigma} & {\bf 0} \\
                   {\bf 0} & {\bf 0} \\
                 \end{array}
               \right)
{\bf V}^{*}\times\\{\bf N}^{-1}{\bf V}\left(
                 \begin{array}{cc}
                   {\bf \Sigma}^{-1} & {\bf 0} \\
                   {\bf 0} & {\bf 0} \\
                 \end{array}
               \right)
{\bf U}^{*}{\bf M}={\rm {\bf X}},\end{multline*}
\begin{multline*}(3M)\,\,({\bf M}{\rm {\bf A}}{\bf X})^{*}  =\left({\bf M}{\bf U}\left(
                 \begin{array}{cc}
                   {\bf \Sigma} & {\bf 0} \\
                   {\bf 0} & {\bf 0} \\
                 \end{array}
               \right)
{\bf V}^{*} {\bf N}^{-1}{\bf V}\left(
                 \begin{array}{cc}
                   {\bf \Sigma}^{-1} & {\bf 0} \\
                   {\bf 0} & {\bf 0} \\
                 \end{array}
               \right)
{\bf U}^{*}{\bf M}\right)^{*}=\\
{\bf M}{\bf U}\left(
                 \begin{array}{cc}
                   {\bf I} & {\bf 0} \\
                   {\bf 0} & {\bf 0} \\
                 \end{array}
               \right)_{m\times m}{\bf U}^{*}{\bf M}={\bf M}{\bf U}\left(
                 \begin{array}{cc}
                   {\bf \Sigma} & {\bf 0} \\
                   {\bf 0} & {\bf 0} \\
                 \end{array}
               \right)
{\bf V}^{*} {\bf N}^{-1}{\bf V}\left(
                 \begin{array}{cc}
                   {\bf \Sigma}^{-1} & {\bf 0} \\
                   {\bf 0} & {\bf 0} \\
                 \end{array}
               \right)
{\bf U}^{*}{\bf M}=\\
{\bf M}{\rm
{\bf A}}{\bf X},\end{multline*}
\begin{multline*}
(4N)\,\,({\bf N}{\bf X}{\rm {\bf A}})^{*}  =\left({\bf N} {\bf N}^{-1}{\bf V}\left(
                 \begin{array}{cc}
                   {\bf \Sigma}^{-1} & {\bf 0} \\
                   {\bf 0} & {\bf 0} \\
                 \end{array}
               \right)
{\bf U}^{*}{\bf M} {\bf U}\left(
                 \begin{array}{cc}
                   {\bf \Sigma} & {\bf 0} \\
                   {\bf 0} & {\bf 0} \\
                 \end{array}
               \right)
{\bf V}^{*}\right)^{*}=\\{\bf N} {\bf N}^{-1}{\bf V}\left(
                 \begin{array}{cc}
                   {\bf I} & {\bf 0} \\
                   {\bf 0} & {\bf 0} \\
                 \end{array}
               \right)_{n\times n}{\bf V}^{*}={\bf N} {\bf N}^{-1}{\bf V}\left(
                 \begin{array}{cc}
                   {\bf \Sigma}^{-1} & {\bf 0} \\
                   {\bf 0} & {\bf 0} \\
                 \end{array}
               \right)
{\bf U}^{*}{\bf M} {\bf U}\left(
                 \begin{array}{cc}
                   {\bf \Sigma} & {\bf 0} \\
                   {\bf 0} & {\bf 0} \\
                 \end{array}
               \right)
{\bf V}^{*}=\\
 {\bf N}{\bf X}{\rm
{\bf A}}.
\end{multline*}$\blacksquare$

\subsection{Limit representations of the weighted Moore-Penrose inverse over the quaternion skew field}
Due to \cite{wei_rep} the following limit representation can be extended to ${\mathbb{H}}$. We  give the proof of the following lemma that different from (\cite{wei_rep}, Corollary 3.4.) and  based on WSVD.
\begin{lemma}\label{lem:lim_rep_weight_mp_A1A}
Let ${\bf A}\in {\mathbb{H}}^{m\times n}_{r}$, and ${\bf M}$ and ${\bf N}$ be positive definite matrices  of order $m$ and $n$, respectively. Then
\begin{equation}\label{eq:lim_rep_weight_mp_A1A}
{\bf A}^{\dag}_{M,N}=\lim_{\lambda \to 0}(\lambda {\bf I}+{\bf A}^{\sharp}{\bf A})^{-1}{\bf A}^{\sharp}.
\end{equation}
where ${\bf A}^{\sharp}={\bf N}^{-1}{\bf A}^{*}{\bf M}$,  $\lambda\in {\mathbb{R}}_{+}$ and ${\mathbb{R}}_{+}$ is the set of all positive real numbers.

\end{lemma}
{\bf Proof}. By Theorems \ref{theor_weig_A} and \ref{tm_weig_A+}, respectively, we have
\[
{\bf A}={\bf U}\begin{pmatrix}
                   {\bf \Sigma}& {\bf 0} \\
                   {\bf 0} & {\bf 0} \\
                 \end{pmatrix}
{\bf V}^{*},
\,\,{\bf A}^{\dag}_{M,N}={\bf N}^{-1}{\bf V}\begin{pmatrix}
                   {\bf \Sigma}^{-1} & {\bf 0} \\
                   {\bf 0} & {\bf 0} \\
                \end{pmatrix}
{\bf U}^{*}{\bf M},
\]
where ${\bf \Sigma}=diag(\sigma_{1}, \sigma_{2},...,\sigma_{r})$, $\sigma_{1}\geq \sigma_{2}\geq...\geq\sigma_{r}>0$ and $\sigma^{2}_{i} \in {\mathbb{R}}$ is the nonzero eigenvalues of ${\bf N}^{-1}{\bf A}^{*}{\bf M}{\bf A}$.     Consider the matrix
\[{\rm {\bf D} }:=\begin{pmatrix}
                   {\bf \Sigma}& {\bf 0} \\
                   {\bf 0} & {\bf 0} \\
                 \end{pmatrix},\]
where ${\rm {\bf D} } = \left( {\sigma _{ij}} \right) \in {\rm
{\mathbb{H}}}_{r}^{m\times n} $ is such that $\sigma _{11} \ge
\sigma _{22} \ge \ldots \ge \sigma _{rr} > \sigma _{r + 1\,r + 1}
= \ldots = \sigma _{qq} = 0$, $q = \min {\left\{ {n,m} \right\}}$. Then
                 \[ {\bf D}^{*} =\begin{pmatrix}
                   {\bf \Sigma}^{*} & {\bf 0} \\
                   {\bf 0} & {\bf 0} \\
                \end{pmatrix},\,\,
 {\bf D}^{+} =\begin{pmatrix}
                   {\bf \Sigma}^{-1} & {\bf 0} \\
                   {\bf 0} & {\bf 0} \\
                \end{pmatrix},
\]
and ${\bf A}={\bf U}{\bf D}{\bf V}^{*}$, ${\bf A}^{\sharp}={\bf N}^{-1}{\bf V}{\bf D}^{*}{\bf U}^{*}{\bf M}$, ${\bf A}^{\dag}_{M,N}={\bf N}^{-1}{\bf V}{\bf D}^{\dag}{\bf U}^{*}{\bf M}$. Since ${\bf N}^{-1}{\bf V}=({\bf V}^{*})^{-1}$, then we have
\begin{multline*}\lambda {\bf I}+{\bf A}^{\sharp}{\bf A}= \lambda {\bf I}+{\bf N}^{-1}{\bf V}{\bf D}^{*}{\bf U}^{*}{\bf M}{\bf U}{\bf D}{\bf V}^{*}=\lambda {\bf I}+({\bf V}^{*})^{-1}{\bf D}^{*}{\bf D}{\bf V}^{*}=\\({\bf V}^{*})^{-1}(\lambda {\bf I}+{\bf D}^{2}){\bf V}^{*}.
\end{multline*}
Further,
\begin{multline*}
(\lambda {\bf I}+{\bf A}^{\sharp}{\bf A})^{-1}{\bf A}^{\sharp}=({\bf V}^{*})^{-1}(\lambda {\bf I}+{\bf D}^{2})^{-1}{\bf V}^{*}{\bf N}^{-1}{\bf V}{\bf D}^{*}{\bf U}^{*}{\bf M}=\\{\bf N}^{-1}{\bf V}(\lambda {\bf I}+{\bf D}^{2})^{-1}{\bf D}^{*}{\bf U}^{*}{\bf M}.
\end{multline*}
 Consider the matrix
\[
\left( {\lambda
{\bf I}+ {\bf D}^{2}} \right)^{ - 1} {\rm {\bf D} }= \left( {{\begin{array}{*{20}c}
 {{\frac{{\sigma _{1}} }{{\sigma _{1}^{2} + \lambda} }}} \hfill & {\ldots}
\hfill & {0} \hfill & {} \hfill & {} \hfill & {} \hfill \\
 {\ldots}  \hfill & {\ldots}  \hfill & {\ldots}  \hfill & {} \hfill & {{\rm
{\bf 0}}} \hfill & {} \hfill \\
 {0} \hfill & {\ldots}  \hfill & {{\frac{{\sigma _{r}} }{{\sigma _{r}^{2}
+ \lambda} }}} \hfill & {} \hfill & { \vdots}  \hfill & {} \hfill \\
 {} \hfill & { \vdots}  \hfill & {} \hfill & { \ddots}  \hfill & {} \hfill &
{} \hfill \\
 {} \hfill & {{\rm {\bf 0}}} \hfill & {} \hfill & {} \hfill & {{\rm {\bf
0}}} \hfill & {} \hfill \\
\end{array}} } \right).
\]
It is obviously that ${\mathop {\lim} \limits_{\lambda \to 0}} \left( {\lambda {\rm {\bf
I}}+{\rm {\bf D} }^{2}} \right)^{ - 1}{\rm
{\bf D}} = {\bf D}^{ \dag} $. Then,
\[\lim_{\lambda \to 0}{\bf N}^{-1}{\bf V}(\lambda {\bf I}+{\bf D}^{2})^{-1}{\bf D}^{*}{\bf U}^{*}{\bf M}={\bf N}^{-1}{\bf V}{\bf D}^{\dag}{\bf U}^{*}{\bf M}={\bf A}^{\dag}_{M,N}.\]
The lemma is proofed.$\blacksquare$

In the following lemma we give another limit representation of ${\bf A}^{+}_{M,N}$.
\begin{lemma}\label{lem:lim_rep_weight_mp_AA1}
Let ${\bf A}\in {\mathbb{H}}^{m\times n}_{r}$, and ${\bf M}$ and ${\bf N}$ be positive definite matrices  of order $m$ and $n$, respectively. Then
\begin{equation}\label{eq:lim_rep_weight_mp_AA1}
{\bf A}^{+}_{M,N}=\lim_{\lambda \to 0}{\bf A}^{\sharp}(\lambda {\bf I}+{\bf A}{\bf A}^{\sharp})^{-1},
\end{equation}
where ${\bf A}^{\sharp}={\bf N}^{-1}{\bf A}^{*}{\bf M}$,  $\lambda\in {\mathbb{R}}_{+}$.

\end{lemma}
{\bf Proof}. The proof is similar to the proof of Lemma \ref{lem:lim_rep_weight_mp_A1A} by using the fact from Theorem \ref{theor_weig_A} that
 the nonzero eigenvalues of ${\bf A}^{\sharp}{\bf A}$ and ${\bf A}{\bf A}^{\sharp}$ coincide.$\blacksquare$
 
It is evidently the following corollary.
\begin{corollary}\label{cor:rep_A1A_ful}  If
${\rm {\bf A}}\in {\mathbb H}^{m\times n} $, then the following
statements are true.
 \begin{itemize}
\item [ i)] If $\rm{rank}\,{\rm {\bf A}} = n$, then ${\rm {\bf A}}_{M,N}^{ +}
= \left( {{\rm {\bf A}}^{\sharp} {\rm {\bf A}}} \right)^{ - 1}{\rm
{\bf A}}^{\sharp}$ .
\item [ ii)] If $\rm{rank}\,{\rm {\bf A}} =
m$, then ${\rm {\bf A}}_{M,N}^{ +}  = {\rm {\bf A}}^{\sharp}\left( {{\rm
{\bf A}}{\rm {\bf A}}^{\sharp} } \right)^{ - 1}.$
\item [ iii)] If $\rm{rank}\,{\rm {\bf A}} = n = m$, then ${\rm {\bf
A}}_{M,N}^{ +}  = {\rm {\bf A}}^{ - 1}$ .
\end{itemize}
\end{corollary}
\section{Determinantal representations of the weighted Moore-Penrose inverse over the quaternion skew field}

Even though  the  eigenvalues of ${\bf A}^{\sharp}{\bf A}$ and ${\bf A}{\bf A}^{\sharp}$ are real and nonnegative, they are not Hermitian in general. Therefor, we consider two cases, when ${\bf A}^{\sharp}{\bf A}$ and ${\bf A}{\bf A}^{\sharp}$ both or one of them are Hermitian, and when they are non-Hermitian.

\subsection{The case of  Hermitian ${\bf A}^{\sharp}{\bf A}$ and ${\bf A}{\bf A}^{\sharp}$.}

Let $({\bf A}^{\sharp}{\bf A}) \in {\mathbb{H}}^{n\times n}$ be Hermitian. It means that $({\bf A}^{\sharp}{\bf A})^{*}=({\bf A}^{\sharp}{\bf A})$. Since ${\bf N}^{-1}$ and ${\bf M}$  are Hermitian, then
 \[({\bf N}^{-1}{\bf A}^{*}{\bf M}{\bf A})^{*}={\bf A}^{*}{\bf M}{\bf A}{\bf N}^{-1}={\bf N}^{-1}{\bf A}^{*}{\bf M}{\bf A}.\]
So, to the matrix $({\bf A}^{\sharp}{\bf A})$ be Hermitian the matrices ${\bf N}^{-1}$ and $({\bf A}^{*}{\bf M}{\bf A})$ should be commutative.
Similarly, to  $({\bf A}{\bf A}^{\sharp})$ be Hermitian the matrices ${\bf M}$ and $({\bf A}{\bf N}^{-1}{\bf A}^{*})$ should be commutative.

Denote by ${\rm {\bf a}}_{.j}^{\sharp} $ and ${\rm {\bf
a}}_{i.}^{\sharp} $ the $j$-th column  and the $i$-th row of  ${\rm
{\bf A}}^{\sharp} $ respectively.
\begin{lemma} \label{kyrc5}
If ${\rm {\bf A}} \in {\rm {\mathbb{H}}}^{m\times n}_r$, then $
 \rank\,\left( {{\rm {\bf A}}^{\sharp} {\rm {\bf A}}}
\right)_{.\,i} \left( {{\rm {\bf a}}_{.j}^{\sharp} }  \right) \le r. $
\end{lemma}
{\bf Proof}. Let's lead  elementary transformations of the
matrix $\left( {{\rm {\bf A}}^{\sharp} {\rm {\bf A}}} \right)_{.\,i}
\left( {{\rm {\bf a}}_{.j}^{\sharp} } \right)$ right-multiplying it by
elementary unimodular matrices ${\rm {\bf P}}_{i\,k} \left( { -
a_{jk}^{}} \right)$, $k \ne j$. The matrix ${\rm {\bf P}}_{\,i\,k}
\left( { - a_{jk}^{}} \right)$ has $-a_{j\,k} $ in the $(i, k)$
entry, 1 in all diagonal entries, and 0 in others. It is the
matrix of an elementary transformation. Right-multiplying  a
matrix by  ${\rm {\bf P}}_{\,i\,k} \left( { - a_{jk}^{}} \right)$,
where $k \ne j$, means adding to $k$-th column its $i$-th column
right-multiplying on $ - a_{jk} $. Then we get

\[
\left( {{\rm {\bf A}}^{\sharp} {\rm {\bf A}}} \right)_{.\,i} \left(
{{\rm {\bf a}}_{.\,j}^{\sharp} }  \right) \cdot {\prod\limits_{k \ne
i} {{\rm {\bf P}}_{i\,k} \left( {-a_{j\,k}}  \right) = {\mathop
{\left( {{\begin{array}{*{20}c}
 {{\sum\limits_{k \ne j} {a_{1k}^{\sharp}  a_{k1}} } } \hfill & {\ldots}  \hfill
& {a_{1j}^{\sharp} }  \hfill & {\ldots}  \hfill & {{\sum\limits_{k \ne
j } {a_{1k}^{\sharp}  a_{kn}}}} \hfill \\
 {\ldots}  \hfill & {\ldots}  \hfill & {\ldots}  \hfill & {\ldots}  \hfill &
{\ldots}  \hfill \\
 {{\sum\limits_{k \ne j} {a_{nk}^{\sharp}  a_{k1}} } } \hfill & {\ldots}  \hfill
& {a_{nj}^{\sharp} }  \hfill & {\ldots}  \hfill & {{\sum\limits_{k \ne
j } {a_{nk}^{\sharp}  a_{kn}}}} \hfill \\
\end{array}} }
\right)}\limits_{i-th}}}}.
\]
The obtained matrix  has the following factorization.
\[
{\mathop {\left( {{\begin{array}{*{20}c}
 {{\sum\limits_{k \ne j} {a_{1k}^{\sharp}  a_{k1}} } } \hfill & {\ldots}  \hfill
& {a_{1j}^{\sharp} }  \hfill & {\ldots}  \hfill & {{\sum\limits_{k \ne
j } {a_{1k}^{\sharp}  a_{kn}} } } \hfill \\
 {\ldots}  \hfill & {\ldots}  \hfill & {\ldots}  \hfill & {\ldots}  \hfill &
{\ldots}  \hfill \\
 {{\sum\limits_{k \ne j} {a_{nk}^{\sharp}  a_{k1}} } } \hfill & {\ldots}  \hfill
& {a_{nj}^{\sharp} }  \hfill & {\ldots}  \hfill & {{\sum\limits_{k \ne
j } {a_{nk}^{\sharp}  a_{kn}} } } \hfill \\
\end{array}} } \right)}\limits_{i-th}}  =
\]
\[
 = \left( {{\begin{array}{*{20}c}
 {a_{11}^{\sharp} }  \hfill & {a_{12}^{\sharp} }  \hfill & {\ldots}  \hfill &
{a_{1m}^{\sharp} }  \hfill \\
 {a_{21}^{\sharp} }  \hfill & {a_{22}^{\sharp} }  \hfill & {\ldots}  \hfill &
{a_{2m}^{\sharp} }  \hfill \\
 {\ldots}  \hfill & {\ldots}  \hfill & {\ldots}  \hfill & {\ldots}  \hfill
\\
 {a_{n1}^{\sharp} }  \hfill & {a_{n2}^{\sharp} }  \hfill & {\ldots}  \hfill &
{a_{nm}^{\sharp} }  \hfill \\
\end{array}} } \right){\mathop {\left( {{\begin{array}{*{20}c}
 {a_{11}}  \hfill & {\ldots}  \hfill & {0} \hfill & {\ldots}  \hfill &
{a_{n1}}  \hfill \\
 {\ldots}  \hfill & {\ldots}  \hfill & {\ldots}  \hfill & {\ldots}  \hfill &
{\ldots}  \hfill \\
 {0} \hfill & {\ldots}  \hfill & {1} \hfill & {\ldots}  \hfill & {0} \hfill
\\
 {\ldots}  \hfill & {\ldots}  \hfill & {\ldots}  \hfill & {\ldots}  \hfill &
{\ldots}  \hfill \\
 {a_{m1}}  \hfill & {\ldots}  \hfill & {0} \hfill & {\ldots}  \hfill &
{a_{mn}}  \hfill \\
\end{array}} } \right)}\limits_{i-th}} j-th.
\]
Denote by ${\rm {\bf \tilde {A}}}: = {\mathop {\left(
{{\begin{array}{*{20}c}
 {a_{11}}  \hfill & {\ldots}  \hfill & {0} \hfill & {\ldots}  \hfill &
{a_{n1}}  \hfill \\
 {\ldots}  \hfill & {\ldots}  \hfill & {\ldots}  \hfill & {\ldots}  \hfill &
{\ldots}  \hfill \\
 {0} \hfill & {\ldots}  \hfill & {1} \hfill & {\ldots}  \hfill & {0} \hfill
\\
 {\ldots}  \hfill & {\ldots}  \hfill & {\ldots}  \hfill & {\ldots}  \hfill &
{\ldots}  \hfill \\
 {a_{m1}}  \hfill & {\ldots}  \hfill & {0} \hfill & {\ldots}  \hfill &
{a_{mn}}  \hfill \\
\end{array}} } \right)}\limits_{i -th}} j -th$.
The matrix ${\rm {\bf \tilde {A}}}$ is obtained from ${\rm {\bf
A}}$ by replacing all entries of the $j$-th row  and of the $i$-th
column with zeroes except that the $(j, i)$ entry equals 1.
Elementary transformations of a matrix do not change its rank and
the rank of a matrix product does not exceed a rank of each
factors. It follows that $\rank\left( {{\rm {\bf A}}^{\sharp} {\rm
{\bf A}}} \right)_{.\,i} \left( {{\rm {\bf a}}_{.j}^{\sharp} } \right)
\le \min \,{\left\{ {\rank{\rm {\bf A}}^{\sharp}, \rank{\rm {\bf
\tilde {A}}}} \right\}}$. It is obviously that $\rank{\rm {\bf
\tilde {A}}} \ge \rank{\rm {\bf A}} = \rank{\rm {\bf A}}^{\sharp} $.
 This completes the
proof.$\blacksquare$

The following lemma has been proved in the same way.
\begin{lemma}
If ${\rm {\bf A}} \in {\rm {\mathbb{H}}}^{m\times n}_r$, then
$\rank\left( {{\rm {\bf A}}{\rm {\bf A}}^{\sharp} } \right)_{.\,i}
\left( {{\rm {\bf a}}_{.j}^{\sharp} } \right) \le r $.
\end{lemma}
Analogues of the characteristic polynomial are considered in the following lemmas.
\begin{lemma}
If ${\rm {\bf A}} \in {\rm {\mathbb{H}}}^{m\times n}$, $t \in
\mathbb{R}$, and $({\bf A}^{\sharp}{\bf A})$ is Hermitian, then

\begin{equation}
\label{kyr4} {\rm{cdet}} _{i} \left( {t{\rm {\bf I}} + {\rm {\bf
A}}^{\sharp} {\rm {\bf A}}} \right)_{.{\kern 1pt} i} \left( {{\rm {\bf
a}}_{.j}^{\sharp} }  \right) = c_{1}^{\left( {ij} \right)} t^{n - 1} +
c_{2}^{\left( {ij} \right)} t^{n - 2} + \ldots + c_{n}^{\left(
{ij} \right)},
\end{equation}

\noindent where $c_{n}^{\left( {ij} \right)} = {\rm{cdet}} _{i}
\left( {{\rm {\bf A}}^{\sharp} {\rm {\bf A}}} \right)_{.\,i} \left(
{{\rm {\bf a}}_{.\,j}^{\sharp} } \right)$  and $c_{k}^{\left( {ij}
\right)} = {\sum\limits_{\beta \in J_{k,\,n} {\left\{ {i}
\right\}}} {{\rm{cdet}} _{i} \left( {\left( {{\rm {\bf A}}^{\sharp}
{\rm {\bf A}}} \right)_{.\,i} \left( {{\rm {\bf a}}_{.\,j}^{\sharp} }
\right)} \right){\kern 1pt}  _{\beta} ^{\beta} } }$ for all $k =
\overline {1,n - 1} $, $i = \overline {1,n}$, and $j = \overline
{1,m}$.
\end{lemma}
{\bf Proof}. Denote by ${\rm
{\bf b}}_{.{\kern 1pt} {\kern 1pt} i} $ the $i$-th column of the Hermitian matrix ${\rm
{\bf A}}^{\sharp} {\rm {\bf A}} = :\left( {b_{ij}}\right)_{n\times n}$.
Consider the Hermitian matrix $\left( {t{\rm {\bf I}} + {\rm {\bf
A}}^{\sharp} {\rm {\bf A}}} \right)_{.{\kern 1pt} {\kern 1pt} i} ({\rm
{\bf b}}_{.{\kern 1pt} {\kern 1pt} i} ) \in {\rm
{\mathbb{H}}}^{n\times n}$. It differs  from $\left( {t{\rm {\bf
I}} + {\rm {\bf A}}^{\sharp} {\rm {\bf A}}} \right)$ an entry $b_{ii}
$. Taking into account Lemma \ref{lem:char_her} we obtain
\begin{equation}
\label{kyr19} \det \left( {t{\rm {\bf I}} + {\rm {\bf A}}^{\sharp}
{\rm {\bf A}}} \right)_{.{\kern 1pt} i} \left( {{\rm {\bf
b}}_{.{\kern 1pt} {\kern 1pt} i}}  \right) = d_{1} t^{n - 1} +
d_{2} t^{n - 2} + \ldots + d_{n},
\end{equation}
where $d_{k} = {\sum\limits_{\beta \in J_{k,\,n} {\left\{ {i}
\right\}}} {\det \left( {{\rm {\bf A}}^{\sharp} {\rm {\bf A}}}
\right){\kern 1pt} {\kern 1pt} _{\beta} ^{\beta} } } $ is the sum
of all principal minors of order $k$ that contain the $i$-th
column for all $k = \overline {1,n - 1} $ and $d_{n} = \det \left(
{{\rm {\bf A}}^{\sharp} {\rm {\bf A}}} \right)$. Therefore, we have
\[{\rm {\bf b}}_{.{\kern 1pt} {\kern 1pt} i} = \left(
{{\begin{array}{*{20}c}
 {{\sum\limits_{l} {a_{1l}^{\sharp}  a_{li}} } } \hfill \\
 {{\sum\limits_{l} {a_{2l}^{\sharp}  a_{li}} } } \hfill \\
 { \vdots}  \hfill \\
 {{\sum\limits_{l} {a_{nl}^{\sharp}  a_{li}} } } \hfill \\
\end{array}} } \right) = {\sum\limits_{l} {{\rm {\bf a}}_{.\,l}^{\sharp}  a_{li}}
},\]
 where ${\rm {\bf a}}_{.{\kern 1pt} {\kern 1pt} l}^{\sharp}  $ is
the $l$th column-vector  of ${\rm {\bf A}}^{\sharp}$ for all
$l=\overline{1,m}$. Taking into account Theorem \ref{theorem:
determinant of hermitian matrix}, Remark
\ref{rem:exp_det} and Proposition \ref{pr:b_into_brak}  we obtain on the one hand
\begin{multline}
\label{kyr20}
  \det \left( {t{\rm {\bf I}} + {\rm {\bf A}}^{\sharp} {\rm {\bf A}}}
\right)_{.{\kern 1pt} i} \left( {{\rm {\bf b}}_{.{\kern 1pt}
{\kern 1pt} i} } \right) = {\rm{cdet}} _{i} \left( {t{\rm {\bf I}}
+ {\rm {\bf A}}^{\sharp} {\rm {\bf A}}} \right)_{.{\kern 1pt} i}
\left( {{\rm {\bf b}}_{.{\kern 1pt} {\kern 1pt} i}}  \right) =\\
   = {\sum\limits_{l} {{\rm{cdet}} _{i} \left( {t{\rm {\bf I}} + {\rm {\bf A}}^{\sharp}{\rm {\bf A}}} \right)_{.{\kern 1pt} l} \left( {{\rm {\bf
a}}_{.{\kern 1pt} {\kern 1pt} l}^{\sharp}  a_{l{\kern 1pt} i}} \right)
= {\sum\limits_{l} {{\rm{cdet}} _{i} \left( {t{\rm {\bf I}} + {\rm
{\bf A}}^{\sharp} {\rm {\bf A}}} \right)_{.{\kern 1pt} i} \left( {{\rm
{\bf a}}_{.{\kern 1pt} {\kern 1pt} l}^{\sharp} }  \right) \cdot {\kern
1pt}} } }} a_{li}
\end{multline}
On the other hand having changed the order of summation,  we get
for all $k = \overline {1,n - 1} $
\begin{multline}
\label{kyr21}
  d_{k} = {\sum\limits_{\beta \in J_{k,\,n} {\left\{ {i} \right\}}}
{\det \left( {{\rm {\bf A}}^{\sharp} {\rm {\bf A}}} \right){\kern 1pt}
{\kern 1pt} _{\beta} ^{\beta} } }  = {\sum\limits_{\beta \in
J_{k,\,n} {\left\{ {i} \right\}}} {{\rm{cdet}} _{i} \left( {{\rm
{\bf A}}^{\sharp} {\rm {\bf A}}} \right){\kern 1pt} {\kern 1pt}
_{\beta} ^{\beta} } }  =\\
  {\sum\limits_{\beta \in J_{k,\,n}
{\left\{ {i} \right\}}} {{\sum\limits_{l} {{\rm{cdet}} _{i} \left(
{\left( {{\rm {\bf A}}^{\sharp} {\rm {\bf A}}{\kern 1pt}}
\right)_{.{\kern 1pt} {\kern 1pt} i} \left( {{\rm {\bf
a}}_{.{\kern 1pt} {\kern 1pt} l}^{\sharp} a_{l\,i}}  \right)}
\right)}} }} {\kern 1pt} _{\beta} ^{\beta}
 =\\ {\sum\limits_{l} { {{\sum\limits_{\beta \in J_{k,\,n} {\left\{ {i}
\right\}}} {{\rm{cdet}} _{i} \left( {\left( {{\rm {\bf A}}^{\sharp}
{\rm {\bf A}}{\kern 1pt}}  \right)_{.{\kern 1pt} i} \left( {{\rm
{\bf a}}_{.{\kern 1pt} {\kern 1pt} l}^{\sharp} }  \right)}
\right){\kern 1pt} _{\beta} ^{\beta} } }} }}  \cdot a_{l{\kern
1pt} i}.
\end{multline}
By substituting (\ref{kyr20}) and (\ref{kyr21}) in (\ref{kyr19}),
and equating factors at $a _ {l \, i} $ when $l = j $, we obtain
the equality (\ref{kyr4}).$\blacksquare$

By analogy can be proved the following lemma.
\begin{lemma}
If ${\rm {\bf A}} \in {\rm {\mathbb{H}}}^{m\times n}$ and $t \in
\mathbb{R}$, and ${\bf A}{\rm {\bf A}}^{\sharp}$ is Hermitian, then
\[ {\rm{rdet}} _{j}  {( t{\rm {\bf I}} + {\rm
{\bf A}}{\rm {\bf A}}^{\sharp} )_{j\,.\,} ({\rm {\bf a}}_{i.}^{\sharp}  )}
   = r_{1}^{\left( {ij} \right)} t^{n
- 1} +r_{2}^{\left( {ij} \right)} t^{n - 2} + \ldots +
r_{n}^{\left( {ij} \right)},
\]

\noindent where  $r_{n}^{\left( {ij} \right)} = {\rm{rdet}} _{j}
{({\rm {\bf A}}{\rm {\bf A}}^{\sharp} )_{j\,.\,} ({\rm {\bf
a}}_{i.\,}^{\sharp} )}$ and $r_{k}^{\left( {ij} \right)} =
{{{\sum\limits_{\alpha \in I_{r,m} {\left\{ {j} \right\}}}
{{\rm{rdet}} _{j} \left( {({\rm {\bf A}}{\rm {\bf A}}^{\sharp}
)_{j\,.\,} ({\rm {\bf a}}_{i.\,}^{\sharp}  )} \right)\,_{\alpha}
^{\alpha} } }}}$ for all $k = \overline {1,n - 1} $, $i =
\overline {1,n}$, and $j = \overline {1,m}$.
\end{lemma}
The following theorem  introduce the determinantal representations of the weighted Moore-Penrose by analogs of the classical adjoint matrix.
\begin{theorem}\label{th:det_rep_A_mn}
Let ${\rm {\bf A}} \in {\rm {\mathbb{H}}}_{r}^{m\times n} $. If  ${\bf A}^{\sharp}{\rm {\bf A}}$ or ${\bf A}{\rm {\bf A}}^{\sharp}$ are Hermitian, then
the weighted Moore-Penrose inverse  ${\rm {\bf A}}_{M,N}^{\dag} = \left( {\tilde{a}_{ij}^{\dag} } \right) \in
{\rm {\mathbb{H}}}_{}^{n\times m} $ possess the following determinantal representations, respectively,
\begin{equation}
\label{kyr5} \tilde{a}_{ij}^{\dag}  = {\frac{{{\sum\limits_{\beta \in
J_{r,\,n} {\left\{ {i} \right\}}} {{\rm{cdet}} _{i} \left( {\left(
{{\rm {\bf A}}^{\sharp} {\rm {\bf A}}} \right)_{\,. \,i} \left( {{\rm
{\bf a}}_{.j}^{\sharp} }  \right)} \right){\kern 1pt} {\kern 1pt}
_{\beta} ^{\beta} } } }}{{{\sum\limits_{\beta \in J_{r,\,\,n}}
{{\left| {\left( {{\rm {\bf A}}^{\sharp} {\rm {\bf A}}} \right){\kern
1pt}  _{\beta} ^{\beta} }  \right|}}} }}},
\end{equation}
or
\begin{equation}
\label{kyr6} \tilde{a}_{ij}^{\dag}  = {\frac{{{\sum\limits_{\alpha \in
I_{r,m} {\left\{ {j} \right\}}} {{\rm{rdet}} _{j} \left( {({\rm
{\bf A}}{\rm {\bf A}}^{\sharp} )_{j\,.\,} ({\rm {\bf a}}_{i.\,}^{\sharp}
)} \right)\,_{\alpha} ^{\alpha} } }}}{{{\sum\limits_{\alpha \in
I_{r,\,m}}  {{\left| {\left( {{\rm {\bf A}}{\rm {\bf A}}^{\sharp} }
\right){\kern 1pt}  _{\alpha} ^{\alpha} } \right|}}} }}}.
\end{equation}
\end{theorem}
{\bf Proof}.  At first we prove (\ref{kyr5}).
By Lemma \ref{lem:lim_rep_weight_mp_A1A},  ${\rm {\bf A}}^{\dag}  = {\mathop
{\lim} \limits_{\alpha \to 0}} \left( {\alpha {\rm {\bf I}} + {\rm
{\bf A}}^{\sharp} {\rm {\bf A}}} \right)^{ - 1}{\rm {\bf A}}^{\sharp }$. Let ${\bf A}^{\sharp}{\rm {\bf A}}$ is Hermitian.
Then matrix $\left( {\alpha {\rm {\bf I}} + {\rm {\bf A}}^{\sharp} {\rm
{\bf A}}} \right) \in {\rm {\mathbb{H}}}^{n\times n}$ is a
full-rank Hermitian matrix. Taking into account Theorem
\ref{kyrc10} it has an inverse, which we represent as a left
inverse matrix
\[
\left( {\alpha {\rm {\bf I}} + {\rm {\bf A}}^{\sharp} {\rm {\bf A}}}
\right)^{ - 1} = {\frac{{1}}{{\det \left( {\alpha {\rm {\bf I}} +
{\rm {\bf A}}^{ * }{\rm {\bf A}}} \right)}}}\left(
{{\begin{array}{*{20}c}
 {L_{11}}  \hfill & {L_{21}}  \hfill & {\ldots}  \hfill & {L_{n1}}  \hfill
\\
 {L_{12}}  \hfill & {L_{22}}  \hfill & {\ldots}  \hfill & {L_{n2}}  \hfill
\\
 {\ldots}  \hfill & {\ldots}  \hfill & {\ldots}  \hfill & {\ldots}  \hfill
\\
 {L_{1n}}  \hfill & {L_{2n}}  \hfill & {\ldots}  \hfill & {L_{nn}}  \hfill
\\
\end{array}} } \right),
\]
\noindent where $L_{ij} $ is a left $ij$-th cofactor  of
$\alpha {\rm {\bf I}} + {\rm {\bf A}}^{\sharp} {\rm {\bf A}}$. Then we
have
\[\begin{array}{l}
   \left( {\alpha {\rm {\bf I}} + {\rm {\bf A}}^{\sharp} {\rm {\bf A}}} \right)^{ -
1}{\rm {\bf A}}^{\sharp}  = \\
  ={\frac{{1}}{{\det \left( {\alpha {\rm
{\bf I}} + {\rm {\bf A}}^{\sharp} {\rm {\bf A}}} \right)}}}\left(
{{\begin{array}{*{20}c}
 {{\sum\limits_{k = 1}^{n} {L_{k1} a_{k1}^{\sharp} } } } \hfill &
{{\sum\limits_{k = 1}^{n} {L_{k1} a_{k2}^{\sharp} } } } \hfill &
{\ldots}
\hfill & {{\sum\limits_{k = 1}^{n} {L_{k1} a_{km}^{\sharp} } } } \hfill \\
 {{\sum\limits_{k = 1}^{n} {L_{k2} a_{k1}^{\sharp} } } } \hfill &
{{\sum\limits_{k = 1}^{n} {L_{k2} a_{k2}^{\sharp} } } } \hfill &
{\ldots}
\hfill & {{\sum\limits_{k = 1}^{n} {L_{k2} a_{km}^{\sharp} } } } \hfill \\
 {\ldots}  \hfill & {\ldots}  \hfill & {\ldots}  \hfill & {\ldots}  \hfill
\\
 {{\sum\limits_{k = 1}^{n} {L_{kn} a_{k1}^{\sharp} } } } \hfill &
{{\sum\limits_{k = 1}^{n} {L_{kn} a_{k2}^{\sharp} } } } \hfill &
{\ldots}
\hfill & {{\sum\limits_{k = 1}^{n} {L_{kn} a_{km}^{\sharp} } } } \hfill \\
\end{array}} } \right).
\end{array}
\]
Using the definition of a left cofactor, we obtain
\begin{equation}
\label{kyr7} {\rm {\bf A}}_{M,N}^{\dag}  = {\mathop {\lim} \limits_{\alpha
\to 0}} \left( {{\begin{array}{*{20}c}
 {{\frac{{{\rm cdet} _{1} \left( {\alpha {\rm {\bf I}} + {\rm {\bf A}}^{\sharp} {\rm
{\bf A}}} \right)_{.1} \left( {{\rm {\bf a}}_{.1}^{\sharp} }
\right)}}{{\det \left( {\alpha {\rm {\bf I}} + {\rm {\bf A}}^{\sharp}
{\rm {\bf A}}} \right)}}}} \hfill & {\ldots}  \hfill &
{{\frac{{{\rm cdet} _{1} \left( {\alpha {\rm {\bf I}} + {\rm {\bf
A}}^{\sharp} {\rm {\bf A}}} \right)_{.1} \left( {{\rm {\bf a}}_{.m}^{\sharp} }  \right)}}{{\det \left( {\alpha {\rm {\bf I}} + {\rm {\bf
A}}^{\sharp} {\rm {\bf A}}} \right)}}}} \hfill \\
 {\ldots}  \hfill & {\ldots}  \hfill & {\ldots}  \hfill \\
 {{\frac{{{\rm cdet} _{n} \left( {\alpha {\rm {\bf I}} + {\rm {\bf A}}^{\sharp} {\rm
{\bf A}}} \right)_{.n} \left( {{\rm {\bf a}}_{.1}^{\sharp} }
\right)}}{{\det \left( {\alpha {\rm {\bf I}} + {\rm {\bf A}}^{\sharp}
{\rm {\bf A}}} \right)}}}} \hfill & {\ldots}  \hfill &
{{\frac{{{\rm cdet} _{n} \left( {\alpha {\rm {\bf I}} + {\rm {\bf
A}}^{\sharp} {\rm {\bf A}}} \right)_{.n} \left( {{\rm {\bf a}}_{.m}^{\sharp} }  \right)}}{{\det \left( {\alpha {\rm {\bf I}} + {\rm {\bf
A}}^{\sharp} {\rm {\bf A}}} \right)}}}} \hfill \\
\end{array}} } \right).
\end{equation}
By Lemma \ref{lem:char_her}, we have $\det \left( {\alpha {\rm {\bf I}}
+ {\rm {\bf A}}^{\sharp} {\rm {\bf A}}} \right) = \alpha ^{n} + d_{1}
\alpha ^{n - 1} + d_{2} \alpha ^{n - 2} + \ldots + d_{n} $, where
$d_{k} = {\sum\limits_{\beta \in J_{k,\,n}}  {{\left| {\left(
{{\rm {\bf A}}^{\sharp} {\rm {\bf A}}} \right){\kern 1pt} {\kern 1pt}
_{\beta} ^{\beta} }  \right|}}} $ is a sum of principal minors of
${\rm {\bf A}}^{\sharp} {\rm {\bf A}}$ of order $k$ for all  $k =
\overline {1,n - 1} $ and $d_{n} = \det {\rm {\bf A}}^{\sharp} {\rm
{\bf A}}$. Since $\rank{\rm {\bf A}}^{\sharp} {\rm {\bf A}} =
\rank{\rm {\bf A}} = r$ and $d_{n} = d_{n - 1} = \ldots = d_{r +
1} = 0$, it follows that
\begin{equation*}\det \left( {\alpha {\rm {\bf I}} + {\rm
{\bf A}}^{\sharp} {\rm {\bf A}}} \right) = \alpha ^{n} + d_{1} \alpha
^{n - 1} + d_{2} \alpha ^{n - 2} + \ldots + d_{r} \alpha ^{n -
r}.
\end{equation*}
By using (\ref{kyr4}),  we have
\begin{equation*}{\rm{cdet}} _{i} \left( {\alpha
{\rm {\bf I}} + {\rm {\bf A}}^{ * }{\rm {\bf A}}} \right)_{.i}
\left( {{\rm {\bf a}}_{.j}^{\sharp} }  \right) = c_{1}^{\left( {ij}
\right)} \alpha ^{n - 1} + c_{2}^{\left( {ij} \right)} \alpha ^{n
- 2} + \ldots + c_{n}^{\left( {ij} \right)}\end{equation*}
 for all $i =
\overline {1,n} $ and $j = \overline {1,m} $, where $c_{k}^{\left(
{ij} \right)} = {\sum\limits_{\beta \in J_{k,\,n} {\left\{ {i}
\right\}}} {{\rm{cdet}} _{i} \left( {({\rm {\bf A}}^{\sharp} {\rm {\bf
A}})_{.\,i} \left( {{\rm {\bf a}}_{.j}^{\sharp} } \right)}
\right){\kern 1pt} {\kern 1pt} _{\beta }^{\beta} } } $ for all $k
= \overline {1,n - 1} $ and $c_{n}^{\left( {ij} \right)} =
{\rm{cdet}} _{i} \left( {{\rm {\bf A}}^{\sharp} {\rm {\bf A}}}
\right)_{.i} \left( {{\rm {\bf a}}_{.j}^{\sharp} }  \right)$.

Now we prove that $c_{k}^{\left( {ij} \right)} = 0$, when $k \ge r
+ 1$ for all $i = \overline {1,n} $, and $j = \overline {1,m} $.
By Lemma \ref{kyrc5} $\rank\left( {{\rm {\bf A}}^{\sharp}{\rm {\bf
A}}} \right)_{.\,i} \left( {{\rm {\bf a}}_{.j}^{\sharp} }  \right) \le
r$, then the matrix $\left( {{\rm {\bf A}}^{\sharp} {\rm {\bf A}}}
\right)_{.\,i} \left( {{\rm {\bf a}}_{.j}^{\sharp} }  \right)$ has no
more $r$ right-linearly independent columns.

Consider $\left( {({\rm {\bf A}}^{\sharp}{\rm {\bf A}})_{\,.\,i}
\left( {{\rm {\bf a}}_{.j}^{\sharp} } \right)} \right){\kern 1pt}
{\kern 1pt} _{\beta} ^{\beta}  $, when $\beta \in J_{k,n} {\left\{
{i} \right\}}$. It is a principal submatrix of  $\left( {{\rm {\bf
A}}^{\sharp} {\rm {\bf A}}} \right)_{.\,i} \left( {{\rm {\bf
a}}_{.j}^{\sharp} } \right)$ of order $k \ge r + 1$. Deleting both its
$i$-th row and column, we obtain a principal submatrix of order $k
- 1$ of  ${\rm {\bf A}}^{\sharp} {\rm {\bf A}}$.
 We denote it by ${\rm {\bf M}}$. The following cases are
possible.

Let $k = r + 1$ and $\det {\rm {\bf M}} \ne 0$. In this case all
columns of $ {\rm {\bf M}} $ are right-linearly independent. The
addition of all of them on one coordinate to columns of
 $\left( {\left( {{\rm {\bf A}}^{\sharp} {\rm {\bf A}}}
\right)_{.\,i} \left( {{\rm {\bf a}}_{.j}^{\sharp} }  \right)} \right)
{\kern 1pt} _{\beta} ^{\beta}$ keeps their right-linear
independence. Hence, they are basis in a matrix $\left(
{\left( {{\rm {\bf A}}^{\sharp} {\rm {\bf A}}} \right)_{\,.\,i} \left(
{{\rm {\bf a}}_{.j}^{\sharp} } \right)} \right){\kern 1pt} {\kern 1pt}
_{\beta} ^{\beta} $, and  the $i$-th column
is the right linear combination of its basic columns. From this by
Theorem \ref{kyrc12}, we get ${\rm{cdet}} _{i} \left( {\left(
{{\rm {\bf A}}^{\sharp} {\rm {\bf A}}} \right)_{\,.\,i} \left( {{\rm
{\bf a}}_{.j}^{\sharp} } \right)} \right){\kern 1pt} {\kern 1pt}
_{\beta} ^{\beta}  = 0$, when $\beta \in J_{k,n} {\left\{ {i}
\right\}}$ and $k \ge r + 1$.

If $k = r + 1$ and $\det {\rm {\bf M}} = 0$, than $p$, ($p < k$),
columns are basis in  ${\rm {\bf M}}$ and in  $\left( {\left(
{{\rm {\bf A}}^{\sharp} {\rm {\bf A}}} \right)_{.\,i} \left( {{\rm
{\bf a}}_{.j}^{\sharp} } \right)} \right){\kern 1pt} {\kern 1pt}
_{\beta} ^{\beta} $. Then due to Theorem
\ref{kyrc12}, we obtain ${\rm{cdet}} _{i} \left( {\left( {{\rm {\bf
A}}^{\sharp} {\rm {\bf A}}} \right)_{\,.\,i} \left( {{\rm {\bf
a}}_{.j}^{\sharp} } \right)} \right){\kern 1pt} {\kern 1pt} _{\beta}
^{\beta}  = 0$ as well.

If $k > r + 1$, then by Theorem \ref{kyrc13} it
follows that $\det {\rm {\bf M}} = 0$ and $p$, ($p < k - 1$),
columns are basis in the both matrices ${\rm {\bf M}}$ and $\left(
{\left( {{\rm {\bf A}}^{\sharp} {\rm {\bf A}}} \right)_{\,.\,i} \left(
{{\rm {\bf a}}_{.j}^{\sharp} }  \right)} \right){\kern 1pt} {\kern
1pt} _{\beta} ^{\beta}  $. Therefore, by Theorem
\ref{kyrc12}, we obtain  ${\rm{cdet}} _{i} \left( {\left(
{{\rm {\bf A}}^{\sharp} {\rm {\bf A}}} \right)_{\,.\,i} \left( {{\rm
{\bf a}}_{.j}^{\sharp} } \right)} \right){\kern 1pt} {\kern 1pt}
_{\beta} ^{\beta}  = 0$.

Thus in all cases, we have ${\rm{cdet}} _{i} \left( {\left( {{\rm
{\bf A}}^{\sharp}{\rm {\bf A}}} \right)_{\,.\,i} \left( {{\rm {\bf
a}}_{.j}^{\sharp} }  \right)} \right){\kern 1pt} {\kern 1pt} _{\beta}
^{\beta}  = 0$, when $\beta \in J_{k,n} {\left\{ {i} \right\}}$
and $r + 1 \le k < n$, and for all $i = \overline {1,n} $
and $j = \overline {1,m} $,
\begin{multline*}c_{k}^{\left( {ij} \right)} = {\sum\limits_{\beta \in J_{k,\,n}
{\left\{ {i} \right\}}} {{\rm{cdet}} _{i} \left( {\left( {{\rm
{\bf A}}^{\sharp} {\rm {\bf A}}} \right)_{\,.\,i} \left( {{\rm {\bf
a}}_{.j}^{\sharp} } \right)} \right) {\kern 1pt} _{\beta} ^{\beta} }
}=0,\\
c_{n}^{\left( {ij} \right)} = {\rm{cdet}} _{i} \left(
{{\rm {\bf A}}^{\sharp} {\rm {\bf A}}} \right)_{.\,i} \left( {{\rm
{\bf a}}_{.j}^{\sharp} }  \right) = 0.
\end{multline*}

Hence, ${\rm{cdet}} _{i} \left( {\alpha {\rm {\bf I}} + {\rm {\bf
A}}^{\sharp} {\rm {\bf A}}} \right)_{.\,i} \left( {{\rm {\bf
a}}_{.\,j}^{\sharp} }  \right) =c_{1}^{\left( {ij} \right)} \alpha ^{n
- 1} + c_{2}^{\left( {ij} \right)} \alpha ^{n - 2} + \ldots +
c_{r}^{\left( {ij} \right)} \alpha ^{n - r}$ for all $i =
\overline {1,n} $ and $j = \overline {1,m} $. By substituting
these values in the matrix from (\ref{kyr7}), we obtain
\[\begin{array}{c}
  {\rm {\bf A}}_{M,N}^{ +}  = {\mathop {\lim} \limits_{\alpha \to 0}} \left(
{{\begin{array}{*{20}c}
 {{\frac{{c_{1}^{\left( {11} \right)} \alpha ^{n - 1} + \ldots +
c_{r}^{\left( {11} \right)} \alpha ^{n - r}}}{{\alpha ^{n} + d_{1}
\alpha ^{n - 1} + \ldots + d_{r} \alpha ^{n - r}}}}} \hfill &
{\ldots}  \hfill & {{\frac{{c_{1}^{\left( {1m} \right)} \alpha ^{n
- 1} + \ldots + c_{r}^{\left( {1m} \right)} \alpha ^{n -
r}}}{{\alpha ^{n} + d_{1} \alpha
^{n - 1} + \ldots + d_{r} \alpha ^{n - r}}}}} \hfill \\
 {\ldots}  \hfill & {\ldots}  \hfill & {\ldots}  \hfill \\
 {{\frac{{c_{1}^{\left( {n1} \right)} \alpha ^{n - 1} + \ldots +
c_{r}^{\left( {n1} \right)} \alpha ^{n - r}}}{{\alpha ^{n} + d_{1}
\alpha ^{n - 1} + \ldots + d_{r} \alpha ^{n - r}}}}} \hfill &
{\ldots}  \hfill & {{\frac{{c_{1}^{\left( {nm} \right)} \alpha ^{n
- 1} + \ldots + c_{r}^{\left( {nm} \right)} \alpha ^{n -
r}}}{{\alpha ^{n} + d_{1} \alpha
^{n - 1} + \ldots + d_{r} \alpha ^{n - r}}}}} \hfill \\
\end{array}} } \right) =\\
  \left( {{\begin{array}{*{20}c}
 {{\frac{{c_{r}^{\left( {11} \right)}} }{{d_{r}} }}} \hfill & {\ldots}
\hfill & {{\frac{{c_{r}^{\left( {1m} \right)}} }{{d_{r}} }}} \hfill \\
 {\ldots}  \hfill & {\ldots}  \hfill & {\ldots}  \hfill \\
 {{\frac{{c_{r}^{\left( {n1} \right)}} }{{d_{r}} }}} \hfill & {\ldots}
\hfill & {{\frac{{c_{r}^{\left( {nm} \right)}} }{{d_{r}} }}} \hfill \\
\end{array}} } \right).
\end{array}
\]

Here $c_{r}^{\left( {ij} \right)} = {\sum\limits_{\beta \in
J_{r,\,n} {\left\{ {i} \right\}}} {{\rm{cdet}} _{i} \left( {\left(
{{\rm {\bf A}}^{\sharp} {\rm {\bf A}}} \right)_{\,.\,i} \left( {{\rm
{\bf a}}_{.j}^{\sharp} }  \right)} \right){\kern 1pt} {\kern 1pt}
_{\beta} ^{\beta} } } $ and $d_{r} = {\sum\limits_{\beta \in
J_{r,\,\,n}} {{\left| {\left( {{\rm {\bf A}}^{\sharp}{\rm {\bf A}}}
\right){\kern 1pt} {\kern 1pt} _{\beta} ^{\beta} } \right|}}} $.
Thus, we have obtained the determinantal representation of ${\rm
{\bf A}}_{M,N}^{ +}  $ by (\ref{kyr5}).
Similarly can be proved the determinantal representation of
${\rm {\bf A}}_{M,N}^{ +} $ by (\ref{kyr6}).$\blacksquare$

\begin{remark}\label{kyrc19}
If $\rank{\rm {\bf A}} = n$, and
$({\bf A}^{\sharp}{\bf A})$ is Hermitian, then by Corollary \ref{cor:rep_A1A_ful}, ${\rm
{\bf A}}_{M,N}^{ +}  = \left( {{\rm {\bf A}}^{\sharp} {\rm {\bf A}}}
\right)^{ - 1}{\rm {\bf A}}^{\sharp} $. Considering $\left( {{\rm {\bf
A}}^{\sharp} {\rm {\bf A}}} \right)^{ - 1}$ as a left inverse, we get
the following representation of  ${\rm {\bf A}}_{M,N}^{ +} $,

\begin{equation}
\label{kyr8}{\rm {\bf A}}_{M,N}^{ \dag}  = {\frac{{1}}{{\rm {ddet} {\rm
{\bf A}}}}}
\begin{pmatrix}
  {{\rm{cdet}} _{1} ({{\rm {\bf A}}^{\sharp} {\rm {\bf
A}}})_{.\,1} \left( {{\rm {\bf a}}_{.\,1}^{\sharp} }  \right)} &
\ldots & {{\rm{cdet}} _{1} ({{\rm {\bf A}}^{\sharp} {\rm {\bf
A}}})_{.\,1} \left(
{{\rm {\bf a}}_{.\,m}^{\sharp} }  \right)} \\
  \ldots & \ldots & \ldots \\
  {{\rm{cdet}} _{n} ({{\rm {\bf A}}^{\sharp} {\rm {\bf A}}})_{.\,n}
 \left( {{\rm {\bf a}}_{.\,1}^{\sharp} }  \right)} & \ldots &
 {{\rm{cdet}} _{n} ({{\rm {\bf A}}^{\sharp} {\rm {\bf
A}}})_{.\,n} \left( {{\rm {\bf a}}_{.\,m}^{\sharp} }  \right)}.
\end{pmatrix}
\end{equation}
If $m > n$, then by Theorem \ref{th:det_rep_A_mn} for ${\rm {\bf A}}_{M,N}^{ +} $
we have (\ref{kyr5}) as well.
\end{remark}
\begin{remark}\label{kyrc17}
If $\rank{\rm {\bf A}} = m$, and
$({\bf A}{\bf A}^{\sharp})$ is Hermitian, then by Corollary \ref{cor:rep_A1A_ful}, ${\rm
{\bf A}}_{M,N}^{\dag}  = {\rm {\bf A}}^{\sharp} \left( {{\rm {\bf A}}{\rm {\bf
A}}^{\sharp} } \right)^{ - 1}$. Considering  $\left(
{{\rm {\bf A}}{\rm {\bf A}}^{\sharp} } \right)^{ - 1}$ as a right
inverse, we get the following representation of ${\rm {\bf A}}_{M,N}^{
+} $,

\begin{equation}
\label{kyr9} {\rm {\bf A}}_{M,N}^{\dag}  = {\frac{{1}}{{\rm {ddet} {{\rm
{\bf A}} }}}}
\begin{pmatrix}
 {{\rm{rdet}} _{1} ({\rm {\bf A}} {\rm {\bf A}}^{\sharp})_{1.} \left( {{\rm {\bf a}}_{1.}^{\sharp} }  \right)} & \ldots &
 {{\rm{rdet}} _{m} ({\rm {\bf A}} {\rm {\bf A}}^{\sharp})_{m.} \left( {{\rm {\bf a}}_{1.}^{\sharp} }  \right)}\\
 \ldots & \ldots & \ldots \\
  {{\rm{rdet}} _{1} ({\rm {\bf A}} {\rm {\bf A}}^{\sharp})_{1.} \left( {{\rm {\bf a}}_{n.}^{\sharp} }  \right)} & \ldots &
  {{\rm{rdet}} _{m} ({\rm {\bf A}} {\rm {\bf A}}^{\sharp})_{m\,.}
\left( {{\rm {\bf a}}_{n\,.}^{\sharp} }  \right)}
\end{pmatrix}.
\end{equation}
If $m < n$, then by Theorem \ref{th:det_rep_A_mn} for ${\rm {\bf A}}_{M,N}^{ +} $
we also have (\ref{kyr6}).

\end{remark}
We obtain determinantal representations of the projection matrices ${\rm {\bf A}}_{M,N}^{ +} {\rm {\bf A}}$ and  ${\rm {\bf A}}{\rm {\bf A}}_{M,N}^{ \dag} $ in the following corollaries.
\begin{corollary}\label{cor:proj_p}
If ${\rm {\bf A}} \in {\rm {\mathbb{H}}}_{r}^{m\times n} $, where
$r < \min {\left\{ {m,n} \right\}}$ or $r = m < n$, and ${\bf A}^{\sharp}{\rm {\bf A}}$ is Hermitian, then the
projection matrix ${\rm {\bf A}}_{M,N}^{\dag} {\rm {\bf A}} = :{\rm {\bf
P}} = \left( {p_{ij}} \right)_{n\times n} $ possess the following
determinantal representation,
\[
p_{ij} = {\frac{{{\sum\limits_{\beta \in J_{r,\,\,n} {\left\{ {i}
\right\}}} {{\rm{cdet}} _{i} \left( {\left( {{\rm {\bf A}}^{\sharp}
{\rm {\bf A}}} \right)_{.\,i} \left({\rm {\bf d}}_{.j} \right)}
\right){\kern 1pt}  _{\beta} ^{\beta} } }}}{{{\sum\limits_{\beta
\in J_{r,\,n}}  {{\left| {\left( {{\rm {\bf A}}^{\sharp} {\rm {\bf
A}}} \right){\kern 1pt}  _{\beta} ^{\beta} } \right|}}} }}},
\]

\noindent where ${\rm {\bf d}}_{.j} $ is the $j$-th column of ${{\rm {\bf
A}}^{\sharp} {\rm {\bf A}}} \in {\rm
{\mathbb{H}}}^{n\times n}$ and for all $i,j = \overline {1,n}$.

\end{corollary}
{\bf Proof}. Representing ${\rm {\bf A}}^{\dag}
$ by (\ref {kyr5}) and right-multiplying  it
  by $ {\rm {\bf A}} $, we obtain the following presentation of  an entry $p _ {ij} $ of $ {\rm {\bf
A}}_{M,N}^{\dag} {\rm {\bf A}} =: {\rm {\bf P}} = \left ({p _ {ij}} \right)
_ {n \times n} $.

\begin{multline*}
  p_{ij} =
{\sum\limits_{k} {{\frac{{{\sum\limits_{\beta \in J_{r,\,n}
{\left\{ {i} \right\}}} {{\rm{cdet}} _{i} \left( {\left( {{\rm
{\bf A}}^{\sharp} {\rm {\bf A}}} \right)_{.\,i} \left({\rm {\bf
a}}_{.\,k}^{\sharp}\right) }  \right) {\kern 1pt} _{\beta} ^{\beta} }
} }}{{{\sum\limits_{\beta \in J_{r,\,n}}  {{\left| {\left( {{\rm
{\bf A}}^{\sharp} {\rm {\bf A}}} \right) {\kern 1pt}
_{\beta} ^{\beta} }  \right|}}} }}}}}  \cdot a_{kj} =\\
   = {\frac{{{\sum\limits_{\beta \in J_{r,\,n} {\left\{ {i} \right\}}}
{{\sum\limits_{k} {{\rm{cdet}} _{i} \left( {\left( {{\rm {\bf
A}}^{\sharp} {\rm {\bf A}}} \right)_{\,.\,i} \left({\rm {\bf
a}}_{.k}^{\sharp}\right) } \right){\kern 1pt} {\kern 1pt} _{\beta}
^{\beta} } }  \cdot \,a_{kj}} } }}{{{\sum\limits_{\beta \in
J_{r,\,\,n}}  {{\left| {\left( {{\rm {\bf A}}^{\sharp} {\rm {\bf A}}}
\right) {\kern 1pt} _{\beta} ^{\beta} }  \right|}}} }}} =
{\frac{{{\sum\limits_{\beta \in J_{r,\,n} {\left\{ {i} \right\}}}
{{\rm{cdet}} _{i} \left( {\left( {{\rm {\bf A}}^{\sharp} {\rm {\bf
A}}} \right)_{.\,i} \left({ {\rm {\bf d}}}_{.\,j}\right)}  \right)
{\kern 1pt} _{\beta} ^{\beta} } }}}{{{\sum\limits_{\beta \in
J_{r,\,n}} {{\left| {\left( {{\rm {\bf A}}^{\sharp} {\rm {\bf A}}}
\right) {\kern 1pt} _{\beta} ^{\beta} } \right|}}} }}},
\end{multline*}
\noindent where ${\bf d}_{.j} $ is the $j$-th column of ${{\rm {\bf
A}}^{\sharp} {\rm {\bf A}}} \in {\rm {\mathbb{H}}}_{}^{n\times n}$ and for all $i, j = \overline {1, n} $.
$\blacksquare$

By analogy can be proved the following corollary.
\begin{corollary}
If ${\rm {\bf A}} \in {\rm {\mathbb{H}}}_{r}^{m\times n}$, where
$r < \min {\left\{ {m,n} \right\}}$ or $r = n < m$, and
$({\bf A}{\bf A}^{\sharp}) \in {\mathbb{H}}^{m\times m}$ is Hermitian, then  the
projection matrix ${\rm {\bf A}}{\rm {\bf A}}_{M,N}^{ \dag} = :{\rm {\bf
Q}} = \left( {q_{ij}} \right)_{m\times m} $ possess the following
determinantal representation,
\[
q_{ij} = {\frac{{{\sum\limits_{\alpha \in I_{r,\,\,m} {\left\{ {j}
\right\}}} {{{\rm{rdet}}_{j} {\left( {({\rm {\bf A}}\,{\rm {\bf A}}^{\sharp}
)_{j{\kern 1pt} .}\, ({\rm {\bf g}}  _{i{\kern 1pt}  .}\, )}
\right) {\kern 1pt} _{\alpha} ^{\alpha} } }}}
}}{{{\sum\limits_{\alpha \in I_{r,\,m}} {{\left| {\left( {{\rm
{\bf A}}{\rm {\bf A}}^{\sharp} } \right){\kern 1pt} _{\alpha
}^{\alpha} }  \right|}}} }}},
\]
\noindent where ${\rm {\bf g}}_{i.} $ is the $i$-th row of $({\rm
{\bf A}} {\rm {\bf A}}^{\sharp})\in {\rm {\mathbb{H}}}^{m\times m}$ and
for all $i,j = \overline {1,m}$.

\end{corollary}

\subsection{The case of  non-Hermitian ${\bf A}^{\sharp}{\bf A}$ and ${\bf A}{\bf A}^{\sharp}$.}
In this subsection we derive  determinantal representations of the weighted Moore-Penrose inverse of  ${\bf A}\in {\mathbb{H}}^{m\times n}$, when
 $({\bf A}{\bf A}^{\sharp}) \in {\mathbb{H}}^{m\times m}$ and $({\bf A}^{\sharp}{\bf A}) \in {\mathbb{H}}^{n\times n}$ are non-Hermitian.

 First, let $({\bf A}^{\sharp}{\bf A}) \in {\mathbb{H}}^{n\times n}$ be non-Hermitian and $\rank({\bf A}^{\sharp}{\bf A})<n$. Consider
$(\lambda {\bf I}+{\bf A}^{\sharp}{\bf A})^{-1}{\bf A}^{\sharp}$.
We have,
\begin{multline}\label{eq:lA1A}
(\lambda {\bf I}+{\bf A}^{\sharp}{\bf A})^{-1}{\bf A}^{\sharp}=(\lambda {\bf I}+{\bf N}^{-1}{\bf A}^{*}{\bf M}{\bf A})^{-1}{\bf A}^{\sharp}= ({\bf N}^{-1}(\lambda {\bf N}+{\bf A}^{*}{\bf M}{\bf A}))^{-1}{\bf A}^{\sharp}=\\
(\lambda {\bf N}+{\bf A}^{*}{\bf M}{\bf A})^{-1}{\bf A}^{*}{\bf M}={\bf N}^{-\frac{1}{2}}(\lambda +{\bf N}^{-\frac{1}{2}}{\bf A}^{*}{\bf M}{\bf A}{\bf N}^{-\frac{1}{2}})^{-1}{\bf N}^{-\frac{1}{2}}{\bf A}^{*}{\bf M}=\\
{\bf N}^{-\frac{1}{2}}\left(\lambda +\left({\bf M}^{\frac{1}{2}}{\bf A}{\bf N}^{-\frac{1}{2}}\right)^{*}{\bf M}^{\frac{1}{2}}{\bf A}{\bf N}^{-\frac{1}{2}}\right)^{-1}\left({\bf N}^{-\frac{1}{2}}{\bf A}^{*}{\bf M}^{\frac{1}{2}}\right){\bf M}^{\frac{1}{2}}
\end{multline}
Since by Lemma \ref{lem:lim_rep_A1}
\[\lim_{\lambda \to 0}\left(\lambda +\left({\bf M}^{\frac{1}{2}}{\bf A}{\bf N}^{-\frac{1}{2}}\right)^{*}{\bf M}^{\frac{1}{2}}{\bf A}{\bf N}^{-\frac{1}{2}}\right)^{-1}\left({\bf N}^{-\frac{1}{2}}{\bf A}^{*}{\bf M}^{\frac{1}{2}}\right)= \left({\bf M}^{\frac{1}{2}}{\bf A}{\bf N}^{-\frac{1}{2}}\right)^{\dag},
\]
then combining (\ref{eq:lim_rep_weight_mp_A1A}) and (\ref{eq:lA1A}), we obtain
\begin{equation}\label{eq:det_rep_A1_mn_ful}
{\rm {\bf A}}_{M,N}^{\dag}= {\bf N}^{-\frac{1}{2}} \left({\bf M}^{\frac{1}{2}}{\bf A}{\bf N}^{-\frac{1}{2}}\right)^{\dag}{\bf M}^{\frac{1}{2}}.
\end{equation}

Denote by $\hat{a}_{ij}$ an $ij$-th entry of  $\hat{{\bf A}}^{\dag}:=\left({\bf M}^{\frac{1}{2}}{\bf A}{\bf N}^{-\frac{1}{2}}\right)^{\dag}$. By determinantal representing (\ref{eq:det_rer_A_c}) for $\hat{{\bf A}}^{\dag}$, we obtain
 \begin{multline*}
\hat{a}^{\dag}_{ij}=   {\frac{{{\sum\limits_{\beta \in
J_{r,\,n} {\left\{ {i} \right\}}} {{\rm{cdet}} _{i} \left( {\left(
{\left({\bf M}^{\frac{1}{2}}{\bf A}{\bf N}^{-\frac{1}{2}}\right)^{ *} \left({\bf M}^{\frac{1}{2}}{\bf A}{\bf N}^{-\frac{1}{2}}\right)} \right)_{\,. \,i}  {\left({\bf m}^{\frac{1}{2}}{\bf a}{\bf n}^{-\frac{1}{2}}\right) }_{.j}^{ *}  } \right){\kern 1pt} {\kern 1pt}
_{\beta} ^{\beta} } } }}{{{\sum\limits_{\beta \in J_{r,\,\,n}}
{{\left| {\left( {\left({\bf M}^{\frac{1}{2}}{\bf A}{\bf N}^{-\frac{1}{2}}\right)^{ *} \left({\bf M}^{\frac{1}{2}}{\bf A}{\bf N}^{-\frac{1}{2}}\right)} \right){\kern
1pt}  _{\beta} ^{\beta} }  \right|}}} }}}=\\
 {\frac{{{\sum\limits_{\beta \in
J_{r,\,n} {\left\{ {i} \right\}}} {{\rm{cdet}} _{i} \left( {\left(
{\left({\bf M}^{\frac{1}{2}}{\bf A}{\bf N}^{-\frac{1}{2}}\right)^{ *} \left({\bf M}^{\frac{1}{2}}{\bf A}{\bf N}^{-\frac{1}{2}}\right)} \right)_{\,. \,i}  {\left({\bf n}^{-\frac{1}{2}}{\bf a}^{ *}{\bf m}^{\frac{1}{2}}\right) }_{.j}  } \right){\kern 1pt} {\kern 1pt}
_{\beta} ^{\beta} } } }}{{{\sum\limits_{\beta \in J_{r,\,\,n}}
{{\left| {\left( {\left({\bf M}^{\frac{1}{2}}{\bf A}{\bf N}^{-\frac{1}{2}}\right)^{ *} \left({\bf M}^{\frac{1}{2}}{\bf A}{\bf N}^{-\frac{1}{2}}\right)} \right){\kern
1pt}  _{\beta} ^{\beta} }  \right|}}} }}}
\end{multline*}
where ${\left({\bf m}^{\frac{1}{2}}{\bf a}{\bf n}^{-\frac{1}{2}}\right) }_{.j}^{ *}  $ denote the $j$-th column of $\left({\bf M}^{\frac{1}{2}}{\bf A}{\bf N}^{-\frac{1}{2}}\right)^{ *}$ for all $j=\overline{1,m}$. By $n^{-\frac{1}{2}}_{ik}$ denote an ${ik}$-th entry of ${\bf N}^{-\frac{1}{2}}$ for all $i,k=\overline{1,n}$, and by $m^{\frac{1}{2}}_{lj}$ denote an ${lj}$-th entry of ${\bf M}^{\frac{1}{2}}$ for all $l,j=\overline{1,m}$, respectively.
Then for the weighted Moore-Penrose inverse  ${\rm {\bf A}}_{M,N}^{ +} = \left( {\tilde{a}_{ij}^{ +} } \right) \in
{\rm {\mathbb{H}}}_{}^{n\times m} $, we have
 \begin{multline}\label{eq:det_rep_non_c}
{\tilde{a}_{ij}^{\dag} }= \sum\limits_{k}^{n}\sum\limits_{l}^{m}n^{-\frac{1}{2}}_{ik}\hat{a}^{\dag}_{kl}m^{\frac{1}{2}}_{lj}=\\
{\frac{\sum\limits_{k}n^{-\frac{1}{2}}_{ik}\cdot {{\sum\limits_{\beta \in
J_{r,\,n} {\left\{ {i} \right\}}} {{\rm{cdet}} _{k} \left( {\left(
{\left({\bf M}^{\frac{1}{2}}{\bf A}{\bf N}^{-\frac{1}{2}}\right)^{ *} \left({\bf M}^{\frac{1}{2}}{\bf A}{\bf N}^{-\frac{1}{2}}\right)} \right)_{\,. \,k}  {\left({\bf n}^{-\frac{1}{2}}{\bf a}^{ *}{\bf m}\right) }_{.j}  } \right){\kern 1pt} {\kern 1pt}
_{\beta} ^{\beta} } } }     }{{{\sum\limits_{\beta \in J_{r,\,\,n}}
{{\left| {\left( {\left({\bf M}^{\frac{1}{2}}{\bf A}{\bf N}^{-\frac{1}{2}}\right)^{ *} \left({\bf M}^{\frac{1}{2}}{\bf A}{\bf N}^{-\frac{1}{2}}\right)} \right){\kern
1pt}  _{\beta} ^{\beta} }  \right|}}} }}}=\\
{\frac{\sum\limits_{k}n^{-\frac{1}{2}}_{ik}{{\sum\limits_{\beta \in
J_{r,\,n} {\left\{ {i} \right\}}} {{\rm{cdet}} _{k} \left( {\left(
{\bf N}^{-\frac{1}{2}}{\bf A}^{ *}{\bf M}{\bf A}{\bf N}^{-\frac{1}{2}}
 \right)_{\,. \,k}  {\left({\bf n}^{-\frac{1}{2}}{\bf a}^{ *}{\bf m}\right) }_{.j}  } \right){\kern 1pt} {\kern 1pt}
_{\beta} ^{\beta} } } }     }{{{\sum\limits_{\beta \in J_{r,\,\,n}}
{{\left| {\left(
{\bf N}^{-\frac{1}{2}}{\bf A}^{ *}{\bf M}{\bf A}{\bf N}^{-\frac{1}{2}}
 \right){\kern
1pt}  _{\beta} ^{\beta} }  \right|}}} }}},
\end{multline}
for all $i=\overline{1,n}$, $j=\overline{1,m}$.

If $\rank({\bf A}^{\sharp}{\bf A})=n$, then by Corollary \ref{cor:rep_A1A_ful}, ${\rm {\bf A}}_{M,N}^{\dag}=({\bf A}^{\sharp}{\bf A})^{-1}{\bf A}^{\sharp}$. So,
\begin{equation}\label{eq:rep_A1_ful1}
{\rm {\bf A}}_{M,N}^{\dag}=\left({\bf N}^{-1}{\bf A}^{*}{\bf M}{\bf A}\right)^{-1}{\bf N}^{-1}{\bf A}^{*}{\bf M}=\left({\bf A}^{*}{\bf M}{\bf A}\right)^{-1}{\bf A}^{*}{\bf M}.
\end{equation}
Since ${\bf A}^{*}{\bf M}{\bf A}$ is Hermitian, then we can use the determinantal representation of a Hermitian inverse matrix (\ref{eq:det_her_inv_L}). Denote ${\bf A}^{*}{\bf M}=:(\hat{a})_{ij}\in
{\rm {\mathbb{H}}}^{n\times m}$. So, we have
\begin{gather}\label{eq:det_rep_non_c_ful}
{\tilde{a}_{ij}^{\dag} }=\frac{\sum_{k=1}^{n}L_{ki}\hat{a}_{kj}}{\det({\bf A}^{*}{\bf M}{\bf A})}=\frac{{\rm{cdet}} _{i}({\bf A}^{*}{\bf M}{\bf A})_{.i}\,\hat{\bf{a}}_{.j}}{\det({\bf A}^{*}{\bf M}{\bf A})}.
\end{gather}
where $\hat{\bf{a}}_{.j}$ is the $j$-th column of ${\bf A}^{*}{\bf M}$ for all $j=\overline{1,m}$.

Now, let $({\bf A}{\bf A}^{\sharp}) \in {\mathbb{H}}^{m\times m}$ be non-Hermitian and $\rank({\bf A}{\bf A}^{\sharp})<m$. By determinantal representing (\ref{eq:det_rer_A_r}) for $\hat{{\bf A}}^{\dag}$, we similarly obtain
\begin{equation*}
\hat{a}^{\dag}_{ij}=
 {\frac{{{\sum\limits_{\alpha \in
I_{r,\,m} {\left\{ {j} \right\}}} {{\rm{rdet}} _{j} \left( {\left(
{\left({\bf M}^{\frac{1}{2}}{\bf A}{\bf N}^{-\frac{1}{2}}\right) \left({\bf M}^{\frac{1}{2}}{\bf A}{\bf N}^{-\frac{1}{2}}\right)} \right)^{ *}_{j. \,}  {\left({\bf n}^{-\frac{1}{2}}{\bf a}^{ *}{\bf m}^{\frac{1}{2}}\right) }_{i.}  } \right){\kern 1pt} {\kern 1pt}
_{\alpha} ^{\alpha} } } }}{{{\sum\limits_{\alpha \in I_{r,\,\,m}}
{{\left| {\left( {\left({\bf M}^{\frac{1}{2}}{\bf A}{\bf N}^{-\frac{1}{2}}\right) \left({\bf M}^{\frac{1}{2}}{\bf A}{\bf N}^{-\frac{1}{2}}\right)^{ *}} \right){\kern
1pt}  _{\alpha} ^{\alpha} }  \right|}}} }}},
\end{equation*}
where ${\left({\bf n}^{-\frac{1}{2}}{\bf a}^{ *}{\bf m}^{\frac{1}{2}}\right) }_{i.}  $ denote the $i$-th row of $\left({\bf N}^{-\frac{1}{2}}{\bf A}^{ *}{\bf M}^{\frac{1}{2}}\right)$ for all $i=\overline{1,n}$.
Finally, we get
 \begin{multline}\label{eq:det_rep_non_r}
{\tilde{a}_{ij}^{\dag} }=\\ \sum\limits_{l}{\frac{{{\sum\limits_{\alpha \in
I_{r,\,m} {\left\{ {l} \right\}}} {{\rm{rdet}} _{l} \left( {\left(
{\left({\bf M}^{\frac{1}{2}}{\bf A}{\bf N}^{-\frac{1}{2}}\right) \left({\bf M}^{\frac{1}{2}}{\bf A}{\bf N}^{-\frac{1}{2}}\right)^{ *}} \right)_{l. \,}  {\left({\bf n}^{-1}{\bf a}^{ *}{\bf m}^{\frac{1}{2}}\right) }_{i.}  } \right)
_{\alpha} ^{\alpha} } } }}{{{\sum\limits_{\alpha \in I_{r,\,\,m}}
{{\left| {\left( {\left({\bf M}^{\frac{1}{2}}{\bf A}{\bf N}^{-\frac{1}{2}}\right) \left({\bf M}^{\frac{1}{2}}{\bf A}{\bf N}^{-\frac{1}{2}}\right)^{ *}} \right){\kern
1pt}  _{\alpha} ^{\alpha} }  \right|}}} }}}\cdot m^{\frac{1}{2}}_{lj}=\\
 {\frac{\sum\limits_{l}{{\sum\limits_{\alpha \in
I_{r,\,m} {\left\{ {l} \right\}}} {{\rm{rdet}} _{l} \left( {\left(
{{\bf M}^{\frac{1}{2}}{\bf A}{\bf N}^{-1}{\bf A}{\bf M}^{\frac{1}{2}}} \right)_{l. \,}  {\left({\bf n}^{-1}{\bf a}^{ *}{\bf m}^{\frac{1}{2}}\right) }_{i.}  } \right)
_{\alpha} ^{\alpha} } } }m^{\frac{1}{2}}_{lj}}{{{\sum\limits_{\alpha \in I_{r,\,\,m}}
{{\left| {\left( {{\bf M}^{\frac{1}{2}}{\bf A}{\bf N}^{-1}{\bf A}{\bf M}^{\frac{1}{2}}} \right){\kern
1pt}  _{\alpha} ^{\alpha} }  \right|}}} }}},
\end{multline}
for all $i=\overline{1,n}$, $j=\overline{1,m}$.

If $\rank({\bf A}{\bf A}^{\sharp})=m$, then by Corollary \ref{cor:rep_A1A_ful}, ${\rm {\bf A}}_{M,N}^{ +}={\bf A}^{\sharp}({\bf A}{\bf A}^{\sharp})^{-1}$. So,
\begin{equation}\label{eq:rep_A1_ful2}
{\rm {\bf A}}_{M,N}^{\dag}={\bf N}^{-1}{\bf A}^{*}{\bf M}\left({\bf A}{\bf N}^{-1}{\bf A}^{*}{\bf M}\right)^{-1}={\bf N}^{-1}{\bf A}^{*}\left({\bf A}{\bf N}^{-1}{\bf A}^{*}\right)^{-1}.
\end{equation}
Since ${\bf A}{\bf N}^{-1}{\bf A}^{*}$ is Hermitian and full-rank, then we can use the determinantal representation of a Hermitian inverse matrix (\ref{eq:det_her_inv_R}). Denote ${\bf N}^{-1}{\bf A}^{*}=:(\check{a})_{ij}\in
{\rm {\mathbb{H}}}^{n\times m}$. So, we have
\begin{gather}\label{eq:det_rep_non_r_ful}
{\tilde{a}_{ij}^{ +} }=\frac{\sum_{k=1}^{n}\hat{a}_{ik}R_{jk}}{\det({\bf A}{\bf N}^{-1}{\bf A}^{*})}=\frac{{\rm{rdet}}_{j}({\bf A}{\bf N}^{-1}{\bf A}^{*})_{j.}\hat{\bf{a}}_{i.}}{\det({\bf A}{\bf N}^{-1}{\bf A}^{*})},
\end{gather}
where $\hat{\bf{a}}_{i.}$ is the $i$-th row of ${\bf N}^{-1}{\bf A}^{*}$ for all $i=\overline{1,n}$.

Thus, we have proved the following theorem.
\begin{theorem}\label{th:det_rep_A_mn}
Let ${\rm {\bf A}} \in {\rm {\mathbb{H}}}_{r}^{m\times n} $. If ${\bf A}^{\sharp}{\rm {\bf A}}$ is non-Hermitian, then
the weighted Moore-Penrose inverse  ${\rm {\bf A}}_{M,N}^{\dag} = \left( {\tilde{a}_{ij}^{\dag} } \right) \in
{\rm {\mathbb{H}}}^{n\times m} $ possess the  determinantal representation (\ref{eq:det_rep_non_c}) if $r<n$, and (\ref{eq:det_rep_non_c_ful}) if $r=n$. If ${\bf A}{\rm {\bf A}}^{\sharp}$ is non-Hermitian , then
  ${\rm {\bf A}}_{M,N}^{\dag} = \left( {\tilde{a}_{ij}^{\dag} } \right) $ possess the  determinantal representation (\ref{eq:det_rep_non_r}) if $r<m$, and (\ref{eq:det_rep_non_r_ful}) if $r=m$.
\end{theorem}
\begin{remark}The equations (\ref{eq:det_rep_A1_mn_ful}), (\ref{eq:rep_A1_ful1}), and (\ref{eq:rep_A1_ful2}) expand the similarly well-known representations \cite{ben} of the weighted Moore-Penrose inverse  to quaternion matrices.
\end{remark}
\section{Determinantal representation of the weighted Moore-Penrose solution of system linear equations}
Consider  a right system linear equation over the quaternion skew field,
\begin{equation}\label{eq:Ax}
{\bf A}{\bf x}={\bf b}
\end{equation}
 where
${\bf A} \in
{\rm {\mathbb{H}}}^{m\times n}$ is the coefficient matrix, ${\bf b} \in
{\rm {\mathbb{H}}}^{m\times 1}$
is a column of constants, and
${\bf x} \in
{\rm {\mathbb{H}}}^{n\times 1}$
is a unknown column.
Due to \cite{song1} we have the following theorem that characterizes the weighted Moore-Penrose solution of (\ref{eq:Ax}).
\begin{theorem}
The right system linear equation (\ref{eq:Ax}) with restriction ${\bf x} \in \mathcal{R}_{r}({\bf A}^{\sharp})$ has the unique solution $ \tilde{{\bf x}}={\rm {\bf A}}_{M,N}^{ +}{\bf b}$.
\end{theorem}
\begin{theorem}\label{th:cr_Ax_her}Let ${\bf A} \in
{\rm {\mathbb{H}}}^{m\times n}$ and ${\bf A}^{\sharp}{\bf A} \in {\mathbb{H}}^{n\times n}$ be Hermitian.
\begin{enumerate}
\item[(i)] If $\rank{\rm {\bf A}} = k \le m < n$, then the weighted Moore-Penrose solution  $\tilde{{\bf x}} = (\tilde{x}_{1} ,\ldots
,\tilde{x}_{n} )^{T}$ of (\ref{eq:Ax}) possess the following determinantal representation
\begin{equation}
\label{eq:cr_Ax_noful} \tilde{x}_{i}= {\frac{{{\sum\limits_{\beta \in
J_{r,\,n} {\left\{ {i} \right\}}} {{\rm{cdet}} _{i} \left( {\left(
{{\rm {\bf A}}^{ \sharp} {\rm {\bf A}}} \right)_{\,.\,i} \left( {{\rm
{\bf f}}} \right)} \right){\kern 1pt} {\kern 1pt} _{\beta}
^{\beta} } } }}{{{\sum\limits_{\beta \in J_{r,\,\,n}} {{\left|
{\left( {{\rm {\bf A}}^{ \sharp} {\rm {\bf A}}} \right){\kern 1pt}
{\kern 1pt} _{\beta} ^{\beta} }  \right|}}} }}},
\end{equation}
 where ${\rm {\bf f}} =  {\bf A}^{ \sharp}\, {\rm {\bf b}},  $ for all $i = \overline {1,n} $.
\item[(ii)]
If $\rank{\rm {\bf A}} = n$, then   for all
$i = \overline {1,n} $ we have
\begin{equation}
\label{eq:cr_Ax_ful} \tilde{x}_{i} = {\frac{{{\rm{cdet}} _{i} \left( {{\rm
{\bf A}}^{\sharp} {\rm {\bf A}}} \right)_{.i} \left( {{\rm {\bf f}}}
\right)}}{{\det {\rm {\bf A}}^{\sharp}{\bf A}}}}.
\end{equation}
\end{enumerate}
\end{theorem}
{\bf Proof}. i) If $\rank{\rm {\bf A}} = k \le m < n$, then by
Theorem  \ref{th:det_rep_A_mn} we can represent ${\rm {\bf A}}_{M,N}^{ +}$ by
(\ref{kyr5}).  By
component-wise of $ \tilde{{\bf x}}={\rm {\bf A}}_{M,N}^{ +}{\bf b}$, we have
\begin{multline*}
  \tilde{x}_{i} =
{\sum\limits_{i=1}^{m} {{\frac{{{\sum\limits_{\beta \in J_{r,\,n}
{\left\{ {i} \right\}}} {{\rm{cdet}} _{i} \left( {\left( {{\rm
{\bf A}}^{\sharp} {\rm {\bf A}}} \right)_{.\,i} \left({\rm {\bf
a}}_{.\,j}^{\sharp}\right) }  \right) {\kern 1pt} _{\beta} ^{\beta} }
} }}{{{\sum\limits_{\beta \in J_{r,\,n}}  {{\left| {\left( {{\rm
{\bf A}}^{\sharp} {\rm {\bf A}}} \right) {\kern 1pt}
_{\beta} ^{\beta} }  \right|}}} }}}}}  \cdot b_{j} =\\
   = {\frac{{{\sum\limits_{\beta \in J_{r,\,n} {\left\{ {i} \right\}}}
{{\sum\limits_{j} {{\rm{cdet}} _{i} \left( {\left( {{\rm {\bf
A}}^{\sharp} {\rm {\bf A}}} \right)_{\,.\,i} \left({\rm {\bf
a}}_{.j}^{\sharp}\right) } \right){\kern 1pt} {\kern 1pt} _{\beta}
^{\beta} } }  \cdot \,b_{j}} } }}{{{\sum\limits_{\beta \in
J_{r,\,\,n}}  {{\left| {\left( {{\rm {\bf A}}^{\sharp} {\rm {\bf A}}}
\right) {\kern 1pt} _{\beta} ^{\beta} }  \right|}}} }}} =
{\frac{{{\sum\limits_{\beta \in J_{r,\,n} {\left\{ {i} \right\}}}
{{\rm{cdet}} _{i} \left( {\left( {{\rm {\bf A}}^{\sharp} {\rm {\bf
A}}} \right)_{.\,i} \left({ {\rm {\bf f}}}\right)}  \right)
{\kern 1pt} _{\beta} ^{\beta} } }}}{{{\sum\limits_{\beta \in
J_{r,\,n}} {{\left| {\left( {{\rm {\bf A}}^{\sharp} {\rm {\bf A}}}
\right) {\kern 1pt} _{\beta} ^{\beta} } \right|}}} }}},
\end{multline*}
\noindent where ${\bf f}= {\rm {\bf A}}^{\sharp}{\bf b}$ and for all $i = \overline {1, n} $.

\noindent ii)
If $\rank{\rm {\bf A}} = n$, then
${\rm {\bf A}}_{M,N}^{ +} $ can be represented by (\ref{kyr8}).  Representing
${\rm {\bf A}}^{ +} {\rm {\bf b}}$ by
component-wise directly gives
(\ref{eq:cr_Ax_ful}).
$\blacksquare$

\begin{theorem}Let ${\bf A} \in
{\rm {\mathbb{H}}}^{m\times n}$ and ${\bf A}^{\sharp}{\bf A} \in {\mathbb{H}}^{n\times n}$ be non-Hermitian.
\begin{enumerate}
\item[(i)] If $\rank{\rm {\bf A}} = k \le m < n$, then the weighted Moore-Penrose solution  $\tilde{{\bf x}} = (\tilde{x}_{1} ,\ldots
,\tilde{x}_{n} )^{T}$ of (\ref{eq:Ax}) possess the following determinantal representation
\begin{equation*}
 \tilde{x}_{i}= {\frac{\sum\limits_{k}n^{-\frac{1}{2}}_{ik}{{\sum\limits_{\beta \in
J_{r,\,n} {\left\{ {i} \right\}}} {{\rm{cdet}} _{k} \left( {\left(
{\bf N}^{-\frac{1}{2}}{\bf A}^{ *}{\bf M}{\bf A}{\bf N}^{-\frac{1}{2}}
 \right)_{\,. \,k}  { \hat{\bf f} }  } \right){\kern 1pt} {\kern 1pt}
_{\beta} ^{\beta} } } }     }{{{\sum\limits_{\beta \in J_{r,\,\,n}}
{{\left| {\left(
{\bf N}^{-\frac{1}{2}}{\bf A}^{ *}{\bf M}{\bf A}{\bf N}^{-\frac{1}{2}}
 \right){\kern
1pt}  _{\beta} ^{\beta} }  \right|}}} }}},
\end{equation*}
 where ${\rm { \hat{\bf f}}} =  \left({\bf N}^{-\frac{1}{2}}{\bf A}^{ *}{\bf M}\right)\, {\rm {\bf b}}$ and  $n^{-\frac{1}{2}}_{ik}$ is an ${ik}$-th entry of ${\bf N}^{-\frac{1}{2}}$ for all $i,k=\overline{1,n}$.
\item[(ii)]
If $\rank{\rm {\bf A}} = n$, then   for all
$i = \overline {1,n} $ we have
\begin{equation*}
 \tilde{x}_{i} =\frac{{\rm{cdet}} _{i}({\bf A}^{*}{\bf M}{\bf A})_{.i}\,{ \check{\bf f}}}{\det({\bf A}^{*}{\bf M}{\bf A})}.
\end{equation*}
where ${ \check{\bf f}}={\bf A}^{*}{\bf M}\,{\bf b}$ for all $j=\overline{1,m}$.
\end{enumerate}
\end{theorem}
{\bf Proof}. The proof is similar to the proof of Theorem \ref{th:cr_Ax_her}  using component-wise representations of ${\rm {\bf A}}_{M,N}^{ +} $ by (\ref{eq:det_rep_non_c}) in the (i) point, and (\ref{eq:det_rep_non_c_ful}) in the (ii) point, respectively.
$\blacksquare$

Consider  a left system linear equation over the quaternion skew field,
\begin{equation}\label{eq:xA}
{\bf x}{\bf A}={\bf b}
\end{equation}
 where
${\bf A} \in
{\rm {\mathbb{H}}}^{m\times n}$ is the coefficient matrix, ${\bf b} \in
{\rm {\mathbb{H}}}^{1\times n}$
is a row of constants, and
${\bf x} \in
{\rm {\mathbb{H}}}^{1\times m}$
is a unknown row.
The following theorem  characterizes the weighted Moore-Penrose solution of (\ref{eq:xA}).
\begin{theorem}
The left system linear equation (\ref{eq:xA}) with restriction ${\bf x} \in \mathcal{R}_{l}({\bf A}^{\sharp})$ has the unique solution $ \tilde{{\bf x}}={\bf b}{\rm {\bf A}}_{M,N}^{ +}$.
\end{theorem}
\begin{theorem}\label{th:cr_xA_her}Let ${\bf A} \in
{\rm {\mathbb{H}}}^{m\times n}$ and ${\bf A}{\bf A}^{\sharp} \in {\mathbb{H}}^{m\times m}$ be Hermitian.
\begin{enumerate}
\item[(i)] If $\rank{\rm {\bf A}} = k \le n < m$, then the weighted Moore-Penrose solution  $\tilde{{\bf x}} = (\tilde{x}_{1} ,\ldots
,\tilde{x}_{m} )$ of (\ref{eq:xA}) possess the following determinantal representation
\begin{equation*}
 \tilde{x}_{j}= {\frac{{{\sum\limits_{\alpha \in
I_{r,\,m} {\left\{ {j} \right\}}} {{\rm{rdet}} _{i} \left( {\left(
{{\rm {\bf A}}^{ \sharp} {\rm {\bf A}}} \right)_{j\,.} \left( {{\rm
{\bf g}}} \right)} \right){\kern 1pt} {\kern 1pt} _{\alpha}
^{\alpha} } } }}{{{\sum\limits_{\alpha \in I_{r,\,\,m}} {{\left|
{\left( {{\rm {\bf A}}^{ \sharp} {\rm {\bf A}}} \right){\kern 1pt}
{\kern 1pt} _{\alpha} ^{\alpha} }  \right|}}} }}},
\end{equation*}
 where ${\rm {\bf g}} =  {\rm {\bf b}} \,{\bf A}^{ \sharp},  $ for all $j = \overline {1,m} $.
\item[(ii)]
If $\rank{\rm {\bf A}} = m$, then   for all
$j = \overline {1,m} $ we have
\begin{equation*}
 \tilde{x}_{j} = {\frac{{{\rm{rdet}} _{j} \left( {{\rm
{\bf A}} {\rm {\bf A}}^{\sharp}} \right)_{j.} \left( {{\rm {\bf g}}}
\right)}}{{\det {\rm {\bf A}}{\bf A}^{\sharp}}}}.
\end{equation*}
\end{enumerate}
\end{theorem}
{\bf Proof}. The proof is similar to the proof of Theorem \ref{th:cr_Ax_her}  using component-wise representations of ${\rm {\bf A}}_{M,N}^{ +} $ by (\ref{kyr6}) in the (i) point, and (\ref{kyr9}) in the (ii) point, respectively.
$\blacksquare$

\begin{theorem}Let ${\bf A} \in
{\rm {\mathbb{H}}}^{m\times n}$ and ${\bf A}{\bf A}^{\sharp} \in {\mathbb{H}}^{m\times m}$ be non-Hermitian.
\begin{enumerate}
\item[(i)] If $\rank{\rm {\bf A}} = k \le n < m$, then the weighted Moore-Penrose solution  $\tilde{{\bf x}} = (\tilde{x}_{1} ,\ldots
,\tilde{x}_{m} )$ of (\ref{eq:xA}) possess the following determinantal representation
\begin{equation*}
 \tilde{x}_{j}= {\frac{\sum\limits_{l}{{\sum\limits_{\alpha \in
I_{r,\,m} {\left\{ {l} \right\}}} {{\rm{rdet}} _{l} \left( {\left(
{{\bf M}^{\frac{1}{2}}{\bf A}{\bf N}^{-1}{\bf A}{\bf M}^{\frac{1}{2}}} \right)_{l. \,}  {\left({ \hat{\bf g}}\right) }  } \right)
_{\alpha} ^{\alpha} } } }m^{\frac{1}{2}}_{lj}}{{{\sum\limits_{\alpha \in I_{r,\,\,m}}
{{\left| {\left( {{\bf M}^{\frac{1}{2}}{\bf A}{\bf N}^{-1}{\bf A}{\bf M}^{\frac{1}{2}}} \right){\kern
1pt}  _{\alpha} ^{\alpha} }  \right|}}} }}},
\end{equation*}
 where ${\rm { \hat{\bf g}}} =  {\rm {\bf b}}\left({\bf N}^{-1}{\bf A}^{ *}{\bf M}^{\frac{1}{2}}\right)$,  $m^{\frac{1}{2}}_{lj}$ is ${lj}$-th entry of ${\bf M}^{\frac{1}{2}}$ for all $l,j=\overline{1,m}$.
\item[(ii)]
If $\rank{\rm {\bf A}} = m$, then   for all
$j = \overline {1,m} $ we have
\begin{equation*}
 \tilde{x}_{j} =\frac{{\rm{rdet}} _{j}({\bf A}{\bf N}^{-1}{\bf A}^{*})_{j.}\,{ \check{\bf g}}}{\det({\bf A}{\bf N}^{-1}{\bf A}^{*})}.
\end{equation*}
where ${ \check{\bf g}}={\bf b}{\bf N}^{-1}{\bf A}^{*}$.
\end{enumerate}
\end{theorem}
{\bf Proof}. The proof is similar to the proof of Theorem \ref{th:cr_Ax_her}  using component-wise representations of ${\rm {\bf A}}_{M,N}^{ +} $ by (\ref{eq:det_rep_non_r}) in the (i) point, and (\ref{eq:det_rep_non_r_ful}) in the (ii) point, respectively.
$\blacksquare$

\section{Examples}

In this section, we give  examples to illustrate our results.

1. Let
us consider the matrices
\begin{gather}\label{eq:A}{\bf A}=\begin{pmatrix}
  1 & i & j \\
  -k & i & 1 \\
 k & j & -i\\
  j & -1 & i
\end{pmatrix},\end{gather} \begin{multline}\label{ex:weights}
{\bf N}^{-1}=\begin{pmatrix}
 23 & 16-2i-2j+10k & -16+10i-2j-2k \\
 16+2i+2j-10k & 29 & -19-i-13j-k \\
-16-10i+2j+2k & -19+i+13j+k & 29\\
\end{pmatrix},\\{\bf M}=\begin{pmatrix}
 2 & k & i & 0 \\
 -k & 2 & 0 & j \\
 -i & 0 & 2 & k\\
 0 & -j & -k & 2
\end{pmatrix}.
\end{multline}
By direct calculation we get that leading principal minors of ${\bf M}$ and ${\bf N}^{-1}$ are all positive. Therefore, due to Proposition \ref{pr:lead_min},  ${\bf M}$ and ${\bf N}^{-1}$ are positive definite matrices. Similarly, by direct calculation of leading principal minors of ${\bf A}^{*}{\bf A}$, we obtain $\rank {\bf A}^{*}{\bf A}=\rank {\bf A}=2$.

 Further,
\begin{gather*}{\bf A}^{\sharp}= {\bf M} {\bf A}^{*}{\bf N}^{-1}=\\\tiny \begin{pmatrix}
 51-12i+25j-24k & -43-18i+39k & -18+26i-30j-38k & 19-i-50j-42k \\
 -32i+17j-37k & -24-50i+26j+24k & -5-24i-56j+k & -38-25i-18j-67k \\
5-6i-50j+11k &44+23i-12j+7k & 30+38i+5j+37k & 18-44i+6j+54k\\
\end{pmatrix}.\end{gather*}
Since,
\begin{gather*} {\bf A}^{\sharp} {\bf A}=\\\small \begin{pmatrix}
 178 & 41+47i+47j+43k & -41+43i+47j+47k \\
41-47i-47j-43k & 176 & -40-46i-42j-46k \\
-41-43i-47j-47k & -40+46i+42j+46k & 176\\
\end{pmatrix}\end{gather*}
are Hermitian, then we shall be obtain  ${\bf A}_{M,N}^{ +}=\left(\tilde{a}_{ij}^{ +}\right)\in {\mathbb{H}}^{3\times 4}$ due to Theorem \ref{th:det_rep_A_mn} by Eq. (\ref{kyr5}).

We have, ${{\sum\limits_{\beta \in J_{2,\,\,3}}
{{\left| {\left( {{\rm {\bf A}}^{\sharp} {\rm {\bf A}}} \right){\kern
1pt}  _{\beta} ^{\beta} }  \right|}}} }=23380+23380+23380=70140,$ and
\begin{gather*}
{{{\sum\limits_{\beta \in
J_{2,\,3} {\left\{ {1} \right\}}} {{\rm{cdet}} _{1} \left( {\left(
{{\rm {\bf A}}^{\sharp} {\rm {\bf A}}} \right)_{\,. \,1} \left( {{\rm
{\bf a}}_{.1}^{\sharp} }  \right)} \right){\kern 1pt} {\kern 1pt}
_{\beta} ^{\beta} } } }}=6680+1670i+3340j-5010k+\\6680-5010i+3340j-1640k=13360-3340i+6680j-6680k.
\end{gather*}
Then, \begin{gather*}\tilde{a}_{11}^{ +}=\frac{8-2i+4j-4k}{42}.\end{gather*}
Similarly, we obtain
\begin{gather*}\begin{array}{ccc}
                 \tilde{a}_{12}^{ +}=\frac{-7-3i+6k}{42}, & \tilde{a}_{13}^{ +}=\frac{-3+4i-5j-6k}{42}, &  \tilde{a}_{14}^{+}=\frac{3-8j-7k}{42},
               \end{array}
\end{gather*}
\begin{gather*}\begin{array}{cccc}
                 \tilde{a}_{21}^{ +}=\frac{-5i+3j-6k}{42}, & \tilde{a}_{22}^{ +}=\frac{-4-8i+2j+2k}{42}, &  \tilde{a}_{23}^{+}=\frac{-1-4i-9j}{42},\tilde{a}_{24}^{+}=\frac{-6-4i-3j-11k}{42},
               \end{array}
\end{gather*}
\begin{gather*}\begin{array}{cccc}
                 \tilde{a}_{31}^{ +}=\frac{-1-i-8j+2k}{42}, & \tilde{a}_{32}^{ +}=\frac{7+4i-2j+k}{42}, &  \tilde{a}_{33}^{+}=\frac{5+6i+j+6k}{42},\tilde{a}_{34}^{+}=\frac{3-7i+j+9k}{42}.
               \end{array}
\end{gather*}
Finally, we obtain
\begin{multline}\label{eq:A_1}{\bf A}_{M,N}^{ +}=\\\frac{1}{42} \small\begin{pmatrix}
8-2i+4j-4k & -7-3i+6k & -3+4i-5j-6k & 3-8j-7k \\
 -5i+3j-6k & -4-8i+2j+2k & -1-4i-9j & -6-4i-3j-11k \\
-1-i-8j+2k &7+4i-2j+k & 5+6i+j+6k & 3-7i+j+9k\\
\end{pmatrix}.\end{multline}

2. Consider the right system of linear equations,
\begin{gather}\label{eq:sys}{\bf A}{\bf x}={\bf b},
\end{gather}
where the coefficient matrix ${\bf A}$ is (\ref{eq:A}) and the column ${\bf b}=(1\,\,0\,\,i\,\,k)^{T}$.
Using (\ref{eq:A_1}), by the matrix method we have for the weighted Moore-Penrose solution ${\bf \tilde{x}}={\bf A}_{M,N}^{ +}{\bf b}$ of (\ref{eq:sys}) with weights ${\bf M}$ and ${\bf N}$ from (\ref{ex:weights}),
\begin{equation}\label{eq:sol}
\tilde{x}_{1}=\frac{11-13i-2j+4k}{42},  \tilde{x}_{2}=\frac{15-9i+7j-3k}{42}, \tilde{x}_{3}=\frac{-16+5i+5j+4k}{42}.
\end{equation}
Now, we shall find the weighted Moore-Penrose solution of (\ref{eq:sys}) by Cramer's rule (\ref{eq:cr_Ax_noful}).
Since
\[{\bf f}= {\rm {\bf A}}^{\sharp}{\bf b}=\begin{pmatrix}
 67-80i-12j+25k  \\
91 -55i+43j-19k \\
-97+30i+31j+24k
\end{pmatrix},\]
then we have
\begin{gather*}
 \tilde{x}_{1}= {\frac{{{\sum\limits_{\beta \in
J_{2,\,3} {\left\{ {i} \right\}}} {{\rm{cdet}} _{i} \left( {\left(
{{\rm {\bf A}}^{ \sharp} {\rm {\bf A}}} \right)_{\,.\,i} \left( {{\rm
{\bf f}}} \right)} \right){\kern 1pt} {\kern 1pt} _{\beta}
^{\beta} } } }}{{{\sum\limits_{\beta \in J_{r,\,\,n}} {{\left|
{\left( {{\rm {\bf A}}^{ \sharp} {\rm {\bf A}}} \right){\kern 1pt}
{\kern 1pt} _{\beta} ^{\beta} }  \right|}}} }}}=\frac{18370-21710i-3340j+6680k}{70140}=\\\frac{11-13i-2j+4k}{42},
\end{gather*}
\begin{gather*}
\tilde{x}_{2}=\frac{25050-15030i+11690j-5010k}{70140}=\frac{15-9i+7j-3k}{42},\\
\tilde{x}_{3}=\frac{-26720+8350i+8350j+6680k}{70140}=\frac{-16+5i+5j+4k}{42}.
\end{gather*}
As we expected, the weighted Moore-Penrose  solutions  by Cramer's rule and be the matrix method coincide.
\section{Conclusion} In this paper, we derive  determinantal representations of the weighted Moore-Penrose by WSVD within the framework of the theory of the noncommutative column-row determinants. Recently, within the framework of the theory of the noncommutative column-row determinants we have been obtain the determinantal representations of the Drazin inverse \cite{ky_qdr} and the weighted Drazin inverse \cite{ky_qwdr} and corresponding determinantal representations of generallized inverse solutions of some matrix equations in \cite{ky_cr_q_1,ky_cr_q_me,ky_dif_eq,ky_min_norm_com,ky_min_norm_qua}.

\end{document}